\documentclass[11pt]{article}
\usepackage{epsfig}
\usepackage{amsmath}
\usepackage{amssymb}
\usepackage{graphicx}
\usepackage[latin1]{inputenc}
\usepackage{color}
\usepackage{epstopdf}

{\catcode `\@=11 \global\let\AddToReset=\@addtoreset}
\AddToReset{equation}{section}

\newtheorem{theorem}{Theorem}[section]
\newtheorem{lemma}{\bf Lemma}[section]

\newtheorem{@definition}{\sc Definition}[section]

\newtheorem{@remark}{\sc Remark}[section]

\newtheorem{@example}{\sc Example}[section]

\newcommand{\beqn}{\begin{displaymath}}
\newcommand{\eeqn}{\end{displaymath}}
\newcommand{\beq}{\begin{equation}}
\newcommand{\eeq}{\end{equation}}

\def\mathsf{\bf}
\def\N{\mathbb{N}}

\def\R{\mathbb{R}}
\def\Z{\mathbb{Z}}

\def\i{\mathrm i}
\def\d{\mathrm d}
\def\e{\mathrm e}

\def\E{\mathrm E}
\def\P{\mathrm P}

\def\text{\mbox}

\def\1{{\bf 1}}

\def\var{\mathrm{Var}} 
\def\cov{\mathrm{Cov}}

\newcommand{\set}[1]{\left\{#1\right\}}

\newcommand{\noi}{\noindent}

\def\limiteloin{\renewcommand{\arraystretch}{0.5}
\begin{array}[t]{c}
\stackrel{{\cal D}}{\longrightarrow} \\
{\scriptstyle n \rightarrow\infty}
\end{array}\renewcommand{\arraystretch}{1}}

\def\limiteloin2{\renewcommand{\arraystretch}{0.5}
\begin{array}[t]{c}
\stackrel{D[0,1]^2}{\longrightarrow} \\
{\scriptstyle n \rightarrow\infty}
\end{array}\renewcommand{\arraystretch}{1}}

\def\limitet{\renewcommand{\arraystretch}{0.5}
\begin{array}[t]{c}
\stackrel{}{\longrightarrow} \\
{\scriptstyle t\rightarrow\infty}
\end{array}\renewcommand{\arraystretch}{1}}

\def\limitet0{\renewcommand{\arraystretch}{0.5}
\begin{array}[t]{c}
\stackrel{}{\longrightarrow} \\
{\scriptstyle t\rightarrow 0}
\end{array}\renewcommand{\arraystretch}{1}}

\newtheorem{thm}{Theorem}[section]
\newtheorem{cor}[thm]{Corollary}

\newtheorem{prop}[thm]{Proposition}
\newtheorem{defn}[thm]{Definition}
\newtheorem{rem}{Remark}[section]

\def\vep{\varepsilon}

\def\nn{\nonumber}

\begin{document}

\title{Contemporaneous aggregation of triangular array \\ of random-coefficient AR(1) processes
}

\author{Anne Philippe$^1$, \  Donata  Puplinskait\.e$^{1,2}$ \ and \ Donatas Surgailis$^2$ }
\date{\today \\  \small
\vskip.2cm
$^1$Universit\'e de Nantes and $^2$Vilnius University }
\maketitle

\begin{abstract}

We discuss contemporaneous aggregation of independent copies of a triangular array of random-coefficient
AR(1) processes  with i.i.d. innovations belonging to the domain of
attraction of an infinitely divisible law $W$. The limiting aggregated process  is shown
to exist, under general assumptions on $W$ and the mixing distribution, and is represented as
a mixed infinitely divisible moving-average  $\{\mathfrak{X}(t) \} $ in (\ref{mix}).
Partial sums process of $\{\mathfrak{X}(t) \} $ is discussed under the assumption $\E W^2 < \infty $ and
a mixing density
regularly varying at the ``unit root'' $x=1$ with exponent $\beta >0$.
We show that the above partial
sums process may exhibit four different limit behaviors depending on $\beta $ and
the L\'evy triplet of $W$.
Finally, we study the disaggregation problem for $\{\mathfrak{X}(t) \} $ in spirit
of Leipus et al. (2006) and obtain the weak consistency of the corresponding
estimator of $\phi(x) $ in a suitable $L_2-$space.

\end{abstract}

\begin{quote}

\end{quote}

%%%%%%%%%%%%%%%%%%%%%%%%%%%%%%%%%%%%%%%%%%%%%%%%%%%%%%%
\section{Introduction}
%%%%%%%%%%%%%%%%%%%%%%%%%%%%%%%%%%%%%%%%%%%%%%%%%%%%%%%

The present paper discusses contemporaneous aggregation of $N$ independent copies
\begin{equation} \label{AR}
X^{(N)}_i(t) = a_i X^{(N)}_i(t-1) + \vep^{(N)}_i(t), \qquad t\in \Z, \quad i=1,2, \cdots, N
\end{equation}
of random-coefficient AR(1) process  $X^{(N)}(t) = a X^{(N)}(t-1) + \vep^{(N)}(t), \ t \in \Z $,
where $\{\vep^{(N)} (t), \, t \in \Z \}, \, N=1,2, \cdots $
is a triangular array of i.i.d. random variables  in the domain of  attraction of an
infinitely divisible law  $W$:
\begin{eqnarray} \label{DID}
\sum_{t=1}^N \vep^{(N)}(t)&\to_{\rm d}&W
\end{eqnarray}
and where $a $ is
a r.v., independent of  $\{\vep^{(N)} (t), t \in \Z \}$ and satisfying
$|a| < 1 $ almost surely (a.s.).
The limit aggregated process $\{\mathfrak{X}(t), t \in \Z \} $ is defined as the limit in distribution:
\begin{equation}
\sum_{i=1}^N X^{(N)}_i(t) \ \to_{\rm fdd} \  \mathfrak{X}(t).\label{limaggre}
\end{equation}
Here and below, $\to_{\rm d}$ and $ \to_{\rm fdd}$ denote the weak convergence of distributions and finite-dimensional
distributions, respectively.
A particular case of (\ref{AR})-(\ref{limaggre}) corresponding to
$\vep^{(N)}(t) = N^{-1/2} \zeta(t), $ where $ \{\zeta(t), t \in \Z \}$ are i.i.d. r.v.'s with
zero mean and finite variance, leads to the classical aggregation scheme
of Robinson (1978), Granger (1980) and a Gaussian limit
process   $\{\mathfrak{X}(t) \}. $  See also
Gon\c calves and Gouri\` eroux (1988), Zaffaroni (2004), Oppenheim and Viano (2004),
Celov et al. (2007), Beran et al. (2010) on aggregation of more general  time series models with finite variance.
Puplinskait\.e and Surgailis (2009, 2010) discussed
aggregation of random-coefficient AR(1) processes with infinite variance and innovations
$\vep^{(N)}(t) = N^{-1/\alpha} \zeta(t), $ where  $ \{\zeta(t), t \in \Z \}$ are i.i.d. r.v.'s
in the domain of attraction
of $\alpha-$stable law $W,\, 0 < \alpha < 2. $  Aggregation and disaggregation of autoregressive random fields
was discussed in Lavancier (2005, 2011), Lavancier et al. (2012),
Puplinskait\.e and Surgailis (2012), Leonenko et al. (2013).

The present paper discusses  the existence and properties
of the limit process $\{\mathfrak{X}(t) \} $ in the general triangular aggregation scheme (\ref{AR})-(\ref{limaggre}). Let us describe our main
results.
Theorem \ref{thmaggre} (Sec.~2)
says that under condition (\ref{DID}) and some mild additional conditions, the limit process in (\ref{limaggre})
exists and is written as a stochastic integral
\begin{equation}
\mathfrak{X}(t) \ := \ \sum_{s \le t} \int_{(-1,1)} x^{t-s} M_s({\d}x), \qquad t \in \Z,
\label{mix}
\end{equation}
where  $\{M_s, s\in \Z\} $ are i.i.d. copies of an infinitely divisible (ID) random measure $M$ on $(-1,1)$ with control measure
$\Phi (\d x) := \P (a \in \d x)$ and L\'evy characteristics $(\mu, \sigma, \pi)$ the same as of r.v. $W$ ($M \sim W$) in (\ref{DID}), i.e.,
for any Borel set $A \subset (-1,1)$
\begin{eqnarray}\label{Mchf}
\E \e^{{\i} \theta M(A) }&=&
\e^{\Phi(A) V(\theta)}, \qquad \theta \in \R.
\end{eqnarray}
Here and in the sequel, $V(\theta)$ denotes the log-characteristic function of r.v. $W$:
\begin{eqnarray} \label{Wchf0}
V(\theta)&:=&\log \E \e^{\i \theta W}  \ = \
 \int_{\R} (\e^{ {\i} \theta y} - 1 - {\i} \theta y \1(|y| \le 1)) \pi (\d y) - \frac{1}{2}\theta^2 \sigma^2 +  \i \theta \mu,
\end{eqnarray}
where $ \mu \in \R, \, \sigma \ge 0$ and $\pi$ is a L\'evy measure (see sec.~2 for details).  In the particular case
when
$W$ is $\alpha-$stable, $0< \alpha \le 2$, Theorem \ref{thmaggre} agrees with Puplinskait\.e and Surgailis (2010, Thm. 2.1). We note that
the process $\{\mathfrak{X}(t) \} $ in (\ref{mix}) is stationary, ergodic and has ID finite-dimensional distributions.
According to
the terminology in
Rajput and Rosinski (1989), (\ref{mix}) is called a {\it mixed ID moving-average}.

Section 3 discusses partial sums limits and long memory properties of the aggregated process  $\{\mathfrak{X}(t) \} $ in (\ref{mix})
under the assumption that the mixing distribution $\Phi$ has a probability density $\phi$  varying regularly at $x=1$ with exponent $\beta >0$:
\begin{equation}\label{bcond0}
\phi(x) \ \sim \ C (1-x)^\beta, \qquad  x \to 1
\end{equation}
for some $C >0$. (\ref{bcond0}) is similar to the assumptions on the mixing distribution in Granger (1980),
Zaffaroni (2004) and other papers. In the finite variance case $\sigma^2_W := {\rm Var}(W) < \infty $ the aggregated process
in (\ref{mix}) is covariance stationary provided $\E (1- a^2)^{-1} < \infty $, with covariance
\begin{equation} \label{covX}
r(t) \ := \ {\rm Cov}(\mathfrak{X}(t),  \mathfrak{X}(0))
\ =\ \sigma^2_W \, \E \big[\sum_{s \le 0} a^{t-s} a^{-s} \big] \ = \
\sigma^2_W \, \E \big[\frac{a^t }{1-a^2} \big], \quad \forall t\in \N
\end{equation}
depending on $\sigma^2_W$ and the mixing distribution only.
 Note also  that the autocorrelation function  of $\mathfrak{X}$ only depends on the law of
$a$.    It is well-known that for $0< \beta <1 $ and $a \in [0,1)$ a.s.,  (\ref{bcond0}) implies that
$r(t) \sim C_1 t^{-\beta} \, (t \to \infty)$ with some $C_1 >0$, in other words, the aggregated process   $\{\mathfrak{X}(t) \} $
has nonsummable covariances $ \sum_{t \in \Z} |r(t)| = \infty, $ or {\it covariance long memory}.

Long memory is often characterized by the limit behavior of  partial sums. According to Cox (1984),
a stationary process
$\{Y_t, t \in \Z \} $ is said to have {\it distributional long memory}
if
there exist some constants $A_n \to \infty \ (n \to
\infty) $ and $B_n $ and a (nontrivial) stochastic  process
$\{J(\tau), \tau \ge 0\}$ with dependent increments
such
that
\begin{equation}
A_n^{-1} \sum_{t=1}^{[n\tau]} (Y_t - B_n) \ \rightarrow_{\rm fdd} \
J(\tau). \label{distLM}
\end{equation}
In the case when $\{J(\tau) \}$ in
(\ref{distLM}) has independent increments, the corresponding process $\{Y_t, t \in \Z \} $ is said to have
{\it distributional short memory}.

The main result of Sec.~3 is Theorem \ref{sums} which shows that under conditions  (\ref{bcond0}) and $\E W^2 < \infty,$  partial
sums of the aggregated   $\{\mathfrak{X}(t) \} $ in (\ref{mix}) may exhibit four different limit behaviors, depending on parameters
$\beta, \sigma $ and the behavior of the L\'evy measure $\pi$ at the origin. Write $W \sim ID_2(\sigma, \pi)$ if
$\E W = 0, \, \E W^2  = \sigma^2 + \int_{\R} x^2 \pi (\d x) < \infty, $ in which case $V(\theta)$ of (\ref{Wchf0}) can be
written as
\begin{eqnarray} \label{Wchf1}
V(\theta)&=&\int_{\R} (\e^{ {\i} \theta y} - 1 - {\i} \theta y ) \pi (\d y) - \frac{1}{2}\theta^2 \sigma^2.
\end{eqnarray}
The  L\'evy measure $\pi$ is completely determined by two nonincreasing functions $\Pi^+ (x) := \pi (\{ u> x\}), \,
\Pi^- (x)  := \pi (\{ u\le  -x\}), \, x>0 $ on $\R_+ = (0,\infty)$.  Assume that there exist $\alpha >0$ and
$c^\pm \ge 0, c^+ + c^- >0$ such that
\begin{equation} \label{limPi2}
\lim_{x \to 0} x^{\alpha } \Pi^+(x) = c^+, \qquad  \lim_{x \to 0} x^{\alpha } \Pi^-(x) = c^-.
\end{equation}
Under these assumptions, the four limit behaviors of $S_n(\tau) := \sum_{t=1}^{[n\tau]} \mathfrak{X}(t)  $ correspond to
the following  parameter regions:
\begin{description}
\item [{\rm (i)}]  \ $0< \beta < 1, \  \sigma >0, $
\item [{\rm (ii)}] \ $0 < \beta < 1, \ \sigma = 0, \ 1+ \beta < \alpha < 2, $
\item [{\rm (iii)}] \ $0 < \beta < 1, \ \sigma = 0, \ 0 < \alpha <  1 + \beta, $
\item [{\rm (iv)}] \ $\beta > 1.$
\end{description}
According to Theorem \ref{sums}, the limit process  of $\{S_n(\tau) \}$,
in the sense of (\ref{distLM}) with $B_n = 0$ and suitably growing $A_n$ in respective cases (i) - (iv) is a
\begin{description}
\item [{\rm (i)}]  \ fractional Brownian motion with parameter $H = 1 - (\beta/2),  $
\item [{\rm (ii)}] \ $\alpha-$stable self-similar process $\Lambda_{\alpha, \beta}$ with dependent increments and
self-similarity parameter $H = 1 - (\beta/\alpha),$ defined in
(\ref{ZZ1}) below,
\item [{\rm (iii)}] \ $(1+\beta)-$stable L\'evy process with independent increments,
\item [{\rm (iv)}] \ Brownian motion.
\end{description}
See Theorem \ref{sums} for  precise formulations.
Accordingly,
the process  $\{\mathfrak{X}(t) \} $ in (\ref{mix}) has distributional long memory in cases (i) and (ii) and distributional
short memory in case (iii). At the same time, $\{\mathfrak{X}(t) \} $   has covariance long memory in all three cases (i)-(iii).
Case (iv) corresponds to distributional and covariance short memory. As $\alpha $ increases from $0$ to $2$, the
L\'evy measure in (\ref{limPi2}) increases its ``mass'' near the origin, the limiting case $\alpha = 2$
corresponding to $\sigma >0$ or a positive ``mass''
at $0$. We see from (i)-(ii) that distributional long memory is related to $\alpha $ being large enough, or
small jumps of the random
measure $M$ having sufficient high intensity.  Note that  the critical exponent $\alpha = 1 + \beta$ separating the
long and short memory ``regimes'' in (ii) and (iii)
decreases  with $\beta, $ which is quite natural since  smaller $\beta $ means the mixing distribution putting more weight
near the unit root $a=1$.

Since aggregation leads to a natural loss of
information about aggregated ``micro'' series, an important statistical problem arises to recover
the lost information from the observed sample of the aggregated  process. In the context
of the AR(1) aggregation scheme (\ref{AR})-(\ref{limaggre}) this leads to the so-called
the disaggregation problem,
or reconstruction of the mixing density $\phi (x) $ from
observed sample $\mathfrak{X}(1), \cdots, \mathfrak{X}(n)$ of the aggregated process in (\ref{mix}).
For Gaussian process (\ref{mix}), the disaggregation problem was investigated in Leipus et al. (2006) and
Celov et al. (2010), who constructed an estimator of the mixing density based on its expansion in
an orthogonal polynomial basis.
In Sec. 4  we extend the results in Leipus et al.~(2006) to the case when the aggregated process
is a mixed ID moving-average of  (\ref{mix}) with finite 4th moment and
obtain
the weak consistency of the mixture density estimator
in a suitable $L_2-$space
(Theorem \ref{disagg2}).

The results of our paper could be developed  in several directions. We expect  that Theorem \ref{sums} can be extended
to the aggregation scheme with common innovations and
to infinite variance ID moving-averages  of (\ref{mix}), generalizing the results in
Puplinskait\.e and Surgailis~(2009, 2010).
An interesting open problem is generalizing Theorem \ref{sums} to the random field set-up
of Lavancier (2010) and Puplinskait\.e and Surgailis~(2012).

In what follows, $C$ stands for a positive constant whose precise value is unimportant and which may change
from line to line.

%%%%%%%%%%%%%%%%%%%%%%%%%%%%%%%%%%%%%%%%%%%%%%%%%%%%%%%
\section{Existence of the limiting aggregated process}
%%%%%%%%%%%%%%%%%%%%%%%%%%%%%%%%%%%%%%%%%%%%%%%%%%%%%%%

Consider random-coefficient AR(1) equation
\begin{equation} \label{AR1}
X(t) = a X(t-1) + \vep(t), \qquad t \in \Z,
\end{equation}
where $\{\vep(t), t \in \Z\}$ are i.i.d. r.v.'s with generic distribution $\vep$, and $a\in (-1,1)$ is a random coefficient independent
of $\{\vep(t), t \in \Z\}$.
The following proposition is easy. See, e.g. Brandt~(1986),
Puplinskait\.e and Surgailis~(2009).

\begin{prop} \label{ARsol} Assume that  $\E |\vep|^p < \infty $ for some $0< p
\le 2$ and $\E \vep = 0 \, (p \ge 1)$. Then there exists a unique strictly stationary solution to the AR(1) equation
(\ref{AR1}) given by the series
\begin{equation}
X(t) \ = \  \sum_{k=0}^\infty a^k  \vep(t-k). \label{AR2}
\end{equation}
The series in (\ref{AR2}) converge conditionally a.s. and in
$L_p$ for any $|a| <1$. Moreover, if
\begin{equation}
\E \Big[\frac{1}{1- |a|} \Big] \ < \  \infty \label{moment}
\end{equation}
then the series in (\ref{AR2}) converges
unconditionally in $L_p$.
\end{prop}

Write $W \sim ID(\mu, \sigma, \pi)$ if r.v. $W$ is infinitely divisible having the log-characteristic function in
(\ref{Wchf0}), where $\mu \in \R, \sigma \ge 0$ and $\pi $ is a measure on $\R$ satisfying $\pi (\{0\}) = 0$ and
$\int_{\R} (x^2 \wedge 1) \pi (\d x)  < \infty $, called the L\'evy measure of $W$. It is well-known
that the distribution of $W$ is completely determined by the (characteristic) triplet $(\mu, \sigma, \pi)$ and vice versa.
See, e.g., Sato (1999).

\begin{defn}  Let $\{ \vep^{(N)}, N \in \N^* \} $ be a sequence of r.v.'s tending to 0 in probability,
 and
$W \sim ID(\mu, \sigma, \pi) $ be  an ID r.v.  We say that the sequence $\{ \vep^{(N)}\}$
belongs to the domain of attraction of $W$, denoted $\{ \vep^{(N)}\} \in D(W), $  if
\begin{equation}
({\cal C}_N(\theta))^N  \ \to  \  \E \e^{\i \theta W}, \qquad \forall \, \theta \in \R, \label{chfconv}
\end{equation}
where ${\cal C}_N(\theta) := \E \exp \{\i \theta \vep^{(N)}\}$, $\theta \in \R $, is the characteristic function of $\vep^{(N)}$.

\end{defn}

\vskip0cm

\begin{rem} {\rm  Sufficient and necessary conditions for  $\{ \vep^{(N)}\} \in D(W) $ in terms of the distribution functions of $\vep^{(N)}$
are well-known. See, e.g., Sato (1999), Feller (1966, vol. 2, Ch.~17). In particular, these conditions include the convergences
\begin{eqnarray}\label{NvepN}
N\P (\vep^{(N)} > x)&\to&\Pi^+(x),  \qquad
N\P (\vep^{(N)} < -x)\ \to \ \Pi^-(x)
\end{eqnarray}
at each continuity point \ $x>0$ \  of \ $\Pi^+$, $\Pi^-$, respectively, where $\Pi^\pm$ are defined as in (\ref{limPi2}).

}
\end{rem}

\begin{rem} {\rm By taking logarithms of both sides, condition (\ref{chfconv}) can be rewritten as
\begin{eqnarray}
N \log {\cal C}_N(\theta)&\to&\log \E \e^{\i \theta W}  \ = \ V(\theta), \quad \forall \, \theta \in \R,
 \label{logchf}
\end{eqnarray}
with the convention that the l.h.s.  of (\ref{logchf}) is defined   for $N> N_0(\theta)$ sufficiently large only, since for a fixed $N$,
the characteristic function ${\cal C}_N(\theta)$ may
vanish at some points $\theta$. In the general case, (\ref{logchf}) can be precised  as follows:
For any $\epsilon >0 $ and any $K >0$ there exists $N_0(K,\epsilon) \in \N^*$ such that
\begin{eqnarray}
\sup_{|\theta| < K} \big|N \log {\cal C}_N(\theta)- V(\theta)\big|&<&\epsilon, \qquad \forall \, N > N_0(K,\epsilon).   \label{Klogchf}
\end{eqnarray}
}
\end{rem}

The following definitions introduce some technical conditions, in addition to  $\{ \vep^{(N)}\} \in D(W) $,  needed  to prove the convergence
towards the aggregated process in (\ref{limaggre}).

\begin{defn}
Let $0< \alpha \le 2$ and $ \{ \vep^{(N)}\}$ be a sequence of r.v.'s. Write $\{ \vep^{(N)}\} \in T(\alpha) $ if
there exists a constant $C$ independent of $N$ and $x$ and such that one of the two following conditions hold:
either
\smallskip

\noindent (i) $\alpha = 2$ and $ \E \vep^{(N)} = 0, \,  N \E (\vep^{(N)})^2 \le C, \ $ or

\smallskip

\noindent (ii) $0< \alpha < 2 $ and $N \P (|\vep^{(N)}| > x) \le Cx^{-\alpha}, \ x > 0; $
moreover, $ \E \vep^{(N)} = 0$ whenever $1 < \alpha < 2, $ while, for $\alpha =1 $ we assume that the
distribution of $\vep^{(N)}$ is symmetric.

\end{defn}

\begin{defn}\label{Tailpi}
Let $0< \alpha \le 2$ and $W \sim ID(\mu, \sigma, \pi)$. Write $W \in {\cal T}(\alpha) $ if
there exists a constant $C$ independent of $x$ and such that one of the two following conditions hold:
either
\smallskip

\noindent (i) $\alpha = 2$ and $ \E W =0, \,  \E W^2 < \infty, $
or

\smallskip

\noindent (ii) $0< \alpha < 2 $ and $\Pi^+(x) + \Pi^-(x) \le Cx^{-\alpha}, \ \forall \, x > 0; $
moreover, $
\E W = 0 $ whenever $1 < \alpha < 2, $ while, for $\alpha =1 $ we assume that the distribution of $W$
is symmetric.

\end{defn}

\begin{cor} Let $\{ \vep^{(N)}\} \in D(W),\, W\sim ID(\mu, \sigma, \pi).$ Assume that
$\{ \vep^{(N)}\} \in T(\alpha) $ for some $0< \alpha \le 2$. Then
$W \in {\cal T}(\alpha).$

\end{cor}

\noi {\it Proof.} Let $\alpha = 2$ and $R_N $ denote the l.h.s. of (\ref{DID}). Then
$R^2_N\to_{\rm d} W^2$ and $\E W^2 \le \liminf_{N \to \infty} \E R^2_N = \liminf_{N \to \infty} N \E (\vep^{(N)})^2
 < \infty $ follows by  Fatou's lemma. Then, relation $\E W = \lim_{N \to \infty} \E R_N  = 0$ follows
by the dominated convergence theorem. For $0< \alpha < 2$, relation $\Pi^\pm (x) \le C x^{-\alpha} $ at each continuity
point $x$ of $\Pi^\pm $ follows from  $\{ \vep^{(N)}\} \in T(\alpha) $ and (\ref{NvepN}) and then extends
to all $x > 0$ by monotonicity. Verification of the remaining properties of $W$ in the cases $1< \alpha < 2 $ and $\alpha =1$
is easy and is omitted.  \hfill $\Box$

\vskip.2cm

The main result of this section is the following theorem. Recall that $\{X_i(t) \equiv X^{(N)}_i(t)\}, \, i=1,2, \cdots, N$ are independent copies of
AR(1) process in (\ref{AR1}) with
i.i.d. innovations $\{\vep(t) \equiv \vep^{(N)}(t)\} $ and  random coefficient $a  \in (-1,1)$. Write $M \sim W$ if $M$ is an ID random
measure on $(-1,1)$ with characteristic function as in (\ref{Mchf})-(\ref{Wchf0}).

\begin{thm} \label{thmaggre} Let condition (\ref{moment}) hold. In addition,
assume
that the generic
sequence
$\{ \vep^{(N)}\}$ belongs to the domain of attraction of ID r.v. $W \sim ID(\mu,\sigma,\pi)$ and there exists an $0< \alpha \le 2$ such that $\{ \vep^{(N)}\} \in T(\alpha).$
Then the limiting aggregated process  $\{\mathfrak{X}(t)\}$ in (\ref{limaggre}) exists. It is stationary, ergodic,
has infinitely divisible finite-dimensional distributions, and a stochastic integral
representation as in (\ref{mix}), where $M \sim W $.

\end{thm}

\noindent {\it Proof.} We follow the proof of Theorem 2.1 in Puplinskait\.e and Surgailis (2010).
Fix $m\ge 1 $ and $\theta(1), \cdots, \theta(m) \in \R$. Denote
$$
\vartheta(s,a) := \sum_{t=1}^m \theta(t) a^{t-s}\1 (s \le t).
$$
Then $\sum_{t=1}^m \theta(t) X^{(N)}_i(t) = \sum_{s \in \Z} \vartheta(s,a_i) \vep^{(N)}_i(s), \, i=1, \cdots, N$, and
\begin{eqnarray}\label{Cconv}
\E \exp\Big\{\i \sum_{i=1}^N \sum_{t=1}^m \theta(t) X^{(N)}_i(t)\Big\}
&=&\Big(\E \exp \Big\{ \i \sum_{t=1}^m \theta(t) X^{(N)}(t) \Big\} \Big)^N
\ =\  \Big(1 + \frac{\Theta(N)}{N}\Big)^N,
\end{eqnarray}
where
\begin{eqnarray*}
\Theta(N)&:=&N\Big(\E \Big[\prod_{s \in \Z} {\cal C}_N (\vartheta(s,a))\Big] -1 \Big).
\end{eqnarray*}
From definitions (\ref{mix}), (\ref{Wchf0}) it follows that
\begin{eqnarray}\label{Xchf}
\E \exp\Big\{\i \sum_{t=1}^m \theta(t) \mathfrak{X}(t)\Big\}
&=&\e^{\Theta}, \quad \text{where} \quad
\Theta\ :=\  \E \sum_{s \in \Z} V(\vartheta(s,a)).
\end{eqnarray}
The convergence in (\ref{limaggre}) to the aggregated process of (\ref{mix}) follows from
(\ref{Cconv}), (\ref{Xchf}) and the limit
\begin{equation} \label{limTh}
\lim_{N \to \infty} \Theta(N) = \Theta,
\end{equation}
which will be proved below.

Note first that  $\sup_{a \in [0,1), s \in \Z} |\vartheta(s,a)| \le \sum_{t=1}^m |\theta(t)| =: K$ is bounded and therefore
the logarithm $\log {\cal C}_N (\vartheta(s,a)) $ is well-defined for $N > N_0(K) $ large enough, see   (\ref{Klogchf}),
and $\Theta(N)$ can be rewritten as
\begin{eqnarray*}
\Theta(N)&=&\E N \Big( \exp \Big\{N^{-1} \sum_{s \in \Z} N \log {\cal C}_N (\vartheta(s,a))\Big\} -1 \Big).
\end{eqnarray*}
Then (\ref{limTh}) follows if we show that %for each $a \in [0,1)$,
\begin{equation} \label{logCconv}
\lim_{N\to \infty} \sum_{s \in \Z} N \log {\cal C}_N (\vartheta(s,a))\ =\ \sum_{s \in \Z} V(\vartheta(s,a)),
\qquad \forall \, a \in (-1,1)
\end{equation}
and
\begin{equation}\label{domin}
\sum_{s \in \Z} \big| N \log {\cal C}_N (\vartheta(s,a)) \big|\ \le \ \frac{C}{1 - |a|^\alpha}, \qquad \forall \, a \in (-1,1),
\end{equation}
where $C$ does not depend on $N, a $.

Let us prove (\ref{domin}). It suffices to check the bound
\begin{eqnarray}
N|1- {\cal C}_N(\theta)| &\le&C|\theta|^\alpha. \label{Cmomtail0}
\end{eqnarray}
Indeed, since  $|{\cal C}_N (\vartheta(s,a)) -1|< \epsilon$ for $N$ large enough (see above), so
$\big| N \log {\cal C}_N (\vartheta(s,a)) \big|\ \le \ CN\big|1 - {\cal C}_N (\vartheta(s,a)) \big|$ and
(\ref{Cmomtail0}) implies
\begin{eqnarray} \label{Vmom}
\sum_{s \in \Z} \big| N \log {\cal C}_N (\vartheta(s,a)) \big|&\le&C\sum_{s \in \Z} |\vartheta (s,a)|^\alpha \ \le \ \frac{C}{1 - |a|^\alpha},
\end{eqnarray}
see Puplinskait\.e and Surgailis (2010, (A.4)), proving  (\ref{domin}).

Consider (\ref{Cmomtail0}) for  $1 < \alpha < 2$. Since $\E \vep^{(N)} = 0$ so
${\cal C}_N(\theta)-1 = \int_{\R} (\e^{\i \theta x} - 1 - \i \theta x) \d F_N(x)$ and
\begin{eqnarray}
N|1- {\cal C}_N(\theta)|
&\le&N\big|\int_{-\infty}^0(\e^{\i \theta x} - 1 - \i \theta x) \d F_N(x)\big|   + N\big|\int_0^{\infty}(\e^{\i \theta x} - 1 - \i \theta x) \d(1- F_N(x))  \big|\nonumber \\
&=&|\theta|\Big(\big|\int_{-\infty}^0 NF_N(x)  (\e^{\i \theta x}- 1) \d x \big|+
\big|\int_0^{\infty} N(1- F_N(x))  (\e^{\i \theta x}- 1) \d x \big|\Big) \nonumber \\
&\le&C|\theta|\int_0^\infty x^{-\alpha} ((|\theta| x) \wedge 1) \d x \
\le \ C|\theta|^\alpha, \label{Cmomtail}
\end{eqnarray}
since $ NF_N(x)\1 (x <0) + N(1-F_N(x))\1 (x>0) \le C|x|^{-\alpha} $ and the integral
$$
\int_0^\infty x^{-\alpha} ((|\theta| x) \wedge 1) \d x = |\theta| \int_0^{1/|\theta|} x^{1-\alpha}  \d x
+ \int_{1/|\theta|}^\infty x^{-\alpha}  \d x  = |\theta|^{\alpha -1} (\frac{1}{2-\alpha} + \frac{1}{\alpha -1})
$$
converges. In the case  $\alpha = 2$, we have
$N|{\cal C}_N(\theta)-1|\le \frac{1}{2} \theta^2 N\E (\vep^{(N)})^2 \le C \theta^2$
and (\ref{Cmomtail0}) follows.

Next, let $0< \alpha < 1$. Then
\begin{eqnarray*}
N|1- {\cal C}_N(\theta)|
&\le&N\int_{-\infty}^0|\e^{\i \theta x} - 1| \d F_N(x)   + N\int_0^{\infty}|\e^{\i \theta x} - 1| \, |\d(1- F_N(x))| \ =: \  I_1 + I_2.
\end{eqnarray*}
Here, $I_1 \le  2N\int_{-\infty}^0   ((|\theta|\, |x|) \wedge 1)  \d F_N(x) =  2N \int_{-\infty}^{-1/|\theta|} \d F_N(x)
+ 2N|\theta| \int_{-1/|\theta|}^0 |x | \d F_N(x)  =: 2(I_{11}+ I_{12}).$  We have
$I_{11} = N F_N(-1/|\theta|) \le C |\theta|^\alpha $ and
\begin{eqnarray*}
I_{12}&=&- |\theta|N \int_{-1/|\theta|}^0 x  \d F_N(x) = - |\theta|N \Big(x F_N(x)\big|^{x=0}_{x = - 1/|\theta|}
- \int_{-1/|\theta|}^0 F_N(x) \d x \Big)\\
&=&|\theta|N\Big(- \frac{F_N(-1/|\theta|)}{|\theta|}  + \int_{-1/|\theta|}^0 F_N(x) \d x \Big)\\
&\le& C|\theta|^\alpha +  C|\theta| \int_{-1/|\theta|}^0 |x|^{-\alpha} \d x \  \le \  C|\theta|^\alpha.
\end{eqnarray*}
Since $I_2$ can be evaluated analogously, this proves (\ref{Cmomtail0}) for $0< \alpha < 1$.

It remains to prove (\ref{Cmomtail0}) for $\alpha = 1$. Since $\int_{\{|x| \le 1/|\theta|\}}x \d F_N(x) = 0$ by symmetry
of $\vep^{(N)}$, so
${\cal C}_N(\theta)-1 = J_1+ J_2+J_3+J_4$, where
$
J_1 := \int_{-\infty}^{-1/|\theta|}(\e^{\i \theta x} - 1) \d F_N(x), \
J_2 := \int_{-1/|\theta|}^0 (\e^{\i \theta x} - 1 - \i \theta x) \d F_N(x), \
J_3 := \int_0^{1/|\theta|} (\e^{\i \theta x} - 1 - \i \theta x) \d F_N(x), \
J_4 := \int_{1/|\theta|}^\infty (\e^{\i \theta x} - 1) \d F_N(x).
$
We have
$N|J_1| \le 2N F_N(-1/|\theta| ) \le C|\theta| $ and a similar bound follows for $J_i, i=2,3,4. $
This proves (\ref{Cmomtail0}).
Then (\ref{logCconv}) and the remaining proof of (\ref{limTh}) and Theorem  \ref{thmaggre}
follow as in Puplinskait\.e and Surgailis (2010, proof of Thm. 2.1).  \hfill $\Box$

\medskip

Theorem \ref{thmaggre} applies in the case of innovations
in the domain of attraction of $\alpha-$stable law, see below.

\begin{defn}
Let $0< \alpha \le 2$ and $\zeta$ be a r.v. Write $\zeta \in D(\alpha) $ if

\smallskip

\noindent (i) $\alpha = 2$ and $\E \zeta = 0, \, \E \zeta^2 < \infty, $ or

\smallskip

\noindent (ii) $0< \alpha < 2 $ and there exist some constants $c_1, c_2 \ge 0, c_1 + c_2 >0$ such that
$$
\lim_{x \to \infty} x^\alpha \P(\zeta > x) = c_1  \quad \text{and} \quad
\lim_{x \to -\infty} |x|^\alpha \P(\zeta \le x) = c_2;
$$
moreover, $ \E \zeta = 0 $ whenever $1 < \alpha < 2, $ while, for $\alpha =1 $ we assume that the
distribution of $\zeta $ is symmetric.
\end{defn}

\begin{cor} Let $\vep^{(N)} = N^{-1/\alpha} \zeta $, where $\zeta \in D(\alpha), \, 0 < \alpha \le 2$. Then
$\{\vep^{(N)}\} \in T(\alpha)$ and $\{\vep^{(N)}\} \in D(W)$, where $W$ is $\alpha-$stable r.v. with the
characteristic function
\begin{equation} \label{Wstable}
\E \e^{{\i}\theta W}\ = \
\e^{-|\theta|^\alpha \omega (\theta; \alpha, c_1, c_2)}, \quad \theta \in \R,
\end{equation}
where
\begin{eqnarray} \label{omega}
\omega(\theta; \alpha, c_1, c_2) &:=& \left\{
\begin{array}{ll} \frac{\Gamma (2-\alpha)}{1-\alpha}\Big((c_{
1} + c_{2}) \cos(\pi \alpha/2)
 - {\i}(c_{1} - c_{2}) {\rm sign}(\theta) \sin (\pi \alpha/2)\Big),  &\mbox{$\alpha \neq 1,2,$} \\
(c_{1} + c_{2}) (\pi/2), &\mbox{$\alpha=1,$} \\
\sigma^2/2, &\mbox{$\alpha=2.$}
\end{array} \right.
\end{eqnarray}
In this case, the statement of Theorem   \ref{thmaggre} coincides with
Puplinskait\.e and Surgailis (2010, Thm. 2.1).

\end{cor}

%%%%%%%%%%%%%%%%%%%%%%%%%%%%%%%%%%%%%%%%%%%%%%%%%%%%%%%
\section{Convergence of the partial sums}
%%%%%%%%%%%%%%%%%%%%%%%%%%%%%%%%%%%%%%%%%%%%%%%%%%%%%%%

In this section we study partial sums limits and distributional long memory property  of the aggregated mixed ID moving-average in (\ref{mix}) under condition (\ref{bcond0}) on the mixing  distribution $\Phi$. More precisely, we shall assume that $\Phi$ has a density
$\phi$ in a vicinity $(1-\epsilon, 1), \, 0< \epsilon <1$ of the unit root such that
\begin{equation} \label{bcond}
\phi(x) \ = \psi(x)\  (1-x)^\beta, \qquad x \in (1-\epsilon, 1),
\end{equation}
where $\beta > 0 $ and $\psi (x) $ is an bounded function having a finite limit $\psi(1) := \lim_{x\to 1} \psi(x) >0$. Notice
that no restrictions on the mixing distribution in the interval $(-1, 1-\epsilon] $ with exception of \eqref{moment} are imposed. We also expect that condition
(\ref{bcond}) can be further relaxed by including a slowly varying factor as $x \to 1$.

%density $\phi$. More precisely, we shall assume that $\phi$ has the form
%\begin{equation} \label{bcond}
%\phi(x) \ = \psi(x)\  (1-x)^\beta, \qquad x \in (0, 1),
%\end{equation}
%where $\beta > 0 $ and $\psi (x) $ is an bounded function having a finite limit $\psi(1) := \lim_{x\to 1} \psi(x) >0$.

Consider an independently scattered $\alpha-$stable random measure $N({\d}x, {\d}s) $ on
$(0,\infty) \times {\R}$ with
control measure
$\nu({\d}x, {\d}s) := \psi(1) x^{\beta -\alpha} {\d}x {\d}s$ and characteristic function
$\E \e^{\i \theta N(A)} =  \e^{-|\theta|^\alpha \omega(\theta; \alpha, c^+,c^-) \nu(A)}, \, \theta \in \R,
$ where
$A\subset (0,\infty) \times {\R}$ is a Borel set with $ \nu (A) < \infty$ and $\omega $ is defined at (\ref{omega}).
For $1< \alpha \le 2, \, 0< \beta < \alpha -1 $, introduce the process
\begin{eqnarray}
\Lambda_{\alpha,\beta}(\tau)&:=&\int_{{\R}_+\times {\R}} \big( f(x,\tau-s)-f(x,-s)\big) N({\d}x, {\d}s), \qquad \tau \ge 0,  \qquad \text{where}
\label{ZZ1} \\
f(x,t)&:=&\begin{cases} 1- {\e}^{-xt}, &\text{if} \  x>0 \ \text{and} \ t>0, \\
0,&\text{otherwise,}
\end{cases} \nonumber
\end{eqnarray}
defined as a stochastic integral with respect to the above random measure $N$.
The process $\Lambda_{\alpha,\beta}$ was introduced in Puplinskait\.e and Surgailis (2010). It
has stationary increments, $\alpha-$stable finite-dimensional distributions,
a.s. continuous sample paths and is self-similar with parameter $H = 1- \frac{\beta}{\alpha} \in (\frac 1 \alpha, 1). $ Note that
for $\alpha = 2$,  $\Lambda_{2,\beta}$ is a fractional Brownian motion. Write
$\rightarrow_{D[0,1]}$ for the weak convergence of random processes in the Skorohod space
$D[0,1]$ endowed with the $J_1-$topology.

\begin{thm} \label{sums} Let $\{\mathfrak{X}(t)\} $ be the aggregated process in (\ref{mix}), where
$M \sim W \sim ID_2(\sigma,\pi)$
and the mixing distribution
satisfies \eqref{bcond}  and \eqref{moment}.

\medskip

\noindent (i) Let $0< \beta < 1$ and  $\sigma >0$.
Then
\begin{equation} \label{iconv}
\frac{1}{n^{1 - \frac{\beta}{2}}} \sum_{t=1}^{[n\tau]} \mathfrak{X}(t) \ \rightarrow_{D[0,1]} \ B_H(\tau),
\end{equation}
where
$B_H$ is a fractional Brownian motion with parameter $H := 1- \frac{\beta}{2}$ and variance
$\E B^2_H(\tau) = \sigma^2 \psi(1) \Gamma (\beta -2) \tau^{2H}$.

\medskip

\noindent (ii) Let $0< \beta < 1, \, \sigma = 0$  and  there exist  $1+ \beta < \alpha < 2$ and  $c^\pm \ge 0,\, c^+ + c^- >0$ such that
(\ref{limPi2}) hold.
Then
\begin{equation} \label{iiconv}
\frac{1}{n^{1 - \frac{\beta}{\alpha}} } \sum_{t=1}^{[n\tau]} \mathfrak{X}(t) \ \rightarrow_{D[0,1]} \ \Lambda_{\alpha, \beta}(\tau),
\end{equation}
where  $\Lambda_{\alpha,\beta}$ is defined in (\ref{ZZ1}).

\medskip

\noindent (iii) Let $0< \beta < 1, \sigma = 0, \, \pi \neq 0$  and there exists $0 < \alpha < 1+ \beta$ such that
\begin{equation}\label{pimom}
\int_{\R} |x|^{\alpha} \pi (\d x) < \infty.
\end{equation}
Then
\begin{equation} \label{iiiconv}
\frac{1}{n^{\frac{1}{1+ \beta}}} \sum_{t=1}^{[n\tau]} \mathfrak{X}(t) \ \rightarrow_{\rm fdd} \ L_{1+\beta}(\tau),
\end{equation}
where   $\{L_{1+ \beta}(\tau), \tau \ge 0\}$ is an $(1+\beta)-$stable L\'evy process with log-characteristic function
given
in (\ref{Jiii}) below.

\medskip

\noindent (iv) Let $\beta> 1$. Then
\begin{equation} \label{ivconv}
\frac{1}{n^{1/2}} \sum_{t=1}^{[n\tau]} \mathfrak{X}(t) \ \rightarrow_{\rm
fdd} \ \sigma_\Phi B(\tau),
\end{equation}
where $B$ is a standard Brownian motion with $\E B^2(1) =1 $ and
$\sigma_\Phi$ is defined in (\ref{sigmaPhi}) below. Moreover, if $\beta >2$ and $\pi$ satisfies (\ref{pimom}) with $\alpha = 4 $, the convergence $\rightarrow_{\rm
fdd}$ in (\ref{ivconv}) can be replaced by $\rightarrow_{D[0,1]}$.

\end{thm}

\begin{rem} {\rm Note that the normalization exponents in Theorem \ref{sums} decrease from (i) to (iv):
\begin{equation}\label{powers}
1 - \frac{\beta}{2}\ > \ 1 - \frac{\beta}{\alpha}\ > \  \frac{1}{1+ \beta}\ > \  \frac{1}{2}.
\end{equation}
Hence, we may conclude that the dependence in the aggregated process decreases
from (i) to (iv).
Also note that while
$\{\mathfrak{X}(t)\} $ has finite variance in all cases (i) - (iv), the limit of its partial
sums may have infinite variance as it happens in (ii) and (iii). Apparently, the finite-dimensional convergence
in (\ref{iiiconv}) cannot be replaced by the convergence in $D[0,1]$ with the $J_1-$topology.
See  Mikosch et al. (2002, p.40),  Leipus and Surgailis (2003, Remark 4.1) for related discussion.

}
\end{rem}

\noindent {\it Proof.} Decompose $\{\mathfrak{X}(t)\} $ in (\ref{mix}) as $\mathfrak{X}(t)= \mathfrak{X}_+(t) + \mathfrak{X}_-(t)$, where
$\mathfrak{X}_+(t) := \sum_{s \le t} \int_{(1-\epsilon,1)} x^{t-s} M_s({\d}x), \, \mathfrak{X}_-(t) :=
\sum_{s \le t} \int_{(-1,1-\epsilon]} x^{t-s} M_s({\d}x)$ and $0<\epsilon <0$ is the same as in (\ref{bcond}).
Let us first show that
\begin{equation}\label{Sminus}
S^-_n :\ = \ \sum_{t=1}^n \mathfrak{X}_-(t)\ = \ O_p(n^{1/2}).
\end{equation}
Using (\ref{covX}), we can write
\begin{eqnarray*}
\E (S^-_n)^2&=&\sigma^2 \E \Big[\sum_{t,s=1}^n  \frac{a^{|t-s|}}{1-a^2} \1(-1< a\le 1-\epsilon) \Big]\
\le\ C\sum_{s=1}^n \E \Big[ \frac{1 - a^{n-s}}{(1-a^2)(1-a)}  \1(-1< a\le 1-\epsilon) \Big] \\
&\le&C(n/\epsilon)\E (1-a^2)^{-1}\ = \  O(n),
\end{eqnarray*}
proving (\ref{Sminus}). We see from (\ref{Sminus}) and (\ref{powers}) that $S^-_n$ is negligible in the proof of (i) - (iii)
since the normalizing constants in these statements grow faster than $n^{1/2}$. Therefore in the subsequent
proofs of finite-dimensional  convergence in (i) - (iii) we can assume w.l.g. that $ \mathfrak{X}(t)= \mathfrak{X}_+(t)$.

\smallskip

\noindent Proof of (i). The statement is true if $\pi = 0,$ or  $W \sim {\cal N}(0, \sigma^2)$. In the case $\pi \ne 0$,
split $\mathfrak{X}(t) = \mathfrak{X}_1(t) + \mathfrak{X}_2(t),$  where
$\mathfrak{X}_1(t), \mathfrak{X}_2(t) $ are defined following the decomposition of the measure $M = M_1 + M_2 $
into independent random measures $M_1 \sim W_1 \sim ID_2(\sigma,0)$ and $M_2 \sim W_2 \sim ID_2(0,\pi)$. Let us prove that
\begin{equation} \label{i0conv}
S_{n2} \ := \ \sum_{t=1}^{n} \mathfrak{X}_2(t) \ = \ o_p(n^{1 - \frac{\beta}{2}}).
\end{equation}
Let $V_2(\theta) := \log \E \e^{\i \theta W_2}  = \int_{\R} (\e^{\i \theta x} - 1 - \i \theta x) \pi (\d x) $. Then
\begin{equation} \label{V0conv}
|V_2(\theta)| \ \le \  C\theta^2 \quad (\forall \, \theta \in \R)\qquad \text{and} \qquad
|V_2(\theta)| \ = \ o(\theta^2) \quad (|\theta| \to \infty).
\end{equation}
Indeed, for any $\delta >0,$  $|V_2(\theta)| \le\ \theta^2 I_1(\delta) + 2 |\theta| I_2(\delta)$,
where $I_1(\delta) := \theta^{-2} \int_{|x| \le \delta} |\e^{\i \theta x} - 1 - \i \theta x| \pi (\d x) \le
\int_{|x| \le \delta} x^2 \pi (\d x) \to 0 \ (\delta \to 0) $ and
$I_2(\delta) := (2|\theta|)^{-1} \int_{|x| > \delta} |\e^{\i \theta x} - 1 - \i \theta x| \pi (\d x) \le
 \int_{|x| > \delta} |x| \pi (\d x) <  \infty \ (\forall \, \delta >0). $ Hence,
(\ref{V0conv}) follows.

Relation (\ref{i0conv}) follows from   $J_n := \log \E \exp\big\{ \i \theta  n^{-1 + \frac{\beta}{2}}   S_{n2}  \big\} = o(1)$.
We have
\begin{eqnarray*}
J_n
&=&\sum_{s\in \Z} \int_0^\epsilon V_2\Big(\theta n^{-1 + \beta/2}  \sum_{t=1}^{n} (1-z)^{t-s} \1(t\ge s)\Big) z^\beta \psi(1-z) \d z \ = \
J_{n1} + J_{n2},
\end{eqnarray*}
where $J_{n1} := \sum_{s\le 0} \int_0^\epsilon V_2(\cdots) z^\beta \psi(1-z) \d z, \, J_{n2} := \sum_{s=1}^n  \int_0^\epsilon
V_2(\cdots) z^\beta \psi(1-z) \d z. $
By change of variables: $n z = w,  n - s +1 = n u $,
$J_{n2}$ can be rewritten as
\begin{eqnarray*}
J_{n2}
&=&\sum_{s=1}^{n} \int_0^\epsilon V_2\Big(\frac{\theta(1 - (1-z)^{n - s +1})}{n^{1 - \beta/2}  z}
\Big) z^\beta \psi(1-z) \d z \\
&=&\frac{1}{n^\beta}\int_{1/n}^{1} \d u \int_0^{\epsilon n}  V_2\Big( \frac{\theta n^{\beta/2} (1-  (1- \frac{w}{n})^{[un]})} {w} \Big)
w^\beta  \psi\Big(1-\frac{w}{n}\Big)\d w
\\
&=&\theta^2 \int_{0}^{1} \d u \int_0^{\infty} G_n(u,w) w^{\beta-2} \psi\Big(1-\frac{w}{n}\Big) \d w,
\end{eqnarray*}
where
\begin{eqnarray*}
G_n(u,w)
&:=&\big(1-  (1- \frac{w}{n})^{[un]}\big)^2 \kappa \Big( \frac{\theta n^{\beta/2} (1-  (1- \frac{w}{n})^{[un]})} {w} \Big) \1(1/n < u < 1,
0< w < \epsilon n)
\end{eqnarray*}
and where $\kappa (\theta) := V_2(\theta)/\theta^2 $ is a bounded function vanishing as $ |\theta| \to \infty $; see (\ref{V0conv}).
Therefore $G_n(u,w) \to 0\, (n \to \infty)$ for any $u \in (0,1], w >0$ fixed. We also have
$|G_n(u,w)| \le C \big(1-  (1- \frac{w}{n})^{[un]}\big)^2 \le C (1 - \e^{-wu})^2=: \bar G(u,w),$
where
$\int_{0}^{1} \d u \int_0^{\infty} \bar G(u,w) w^{\beta-2} \d w < \infty. $ Thus, $J_{n2} = o(1)$ follows
by the dominated convergence theorem.  The proof $J_{n1} = o(1)$ using (\ref{V0conv}) follows  by a similar argument.
This proves $J_n = o(1),$ or (\ref{i0conv}).
The tightness of the partial sums process in $D[0,1]$ follows from $\beta < 1 $ and
Kolmogorov's criterion since $ \E \big(\sum_{t=1}^{n} \mathfrak{X}(t)\big)^2 = O(n^{2-\beta}), $  the last relation
being an easy consequence of $r(t) = O(t^{-\beta}),$ see
(\ref{covX}) and the discussion below it.

\medskip

\noi Proof of (ii). Let $S_n(\tau) := \sum_{t=1}^{[n\tau]} \mathfrak{X}(t)$.
Let us prove that for any $0< \tau_1 < \cdots < \tau_m \le 1, \, \theta_1\in \R, \cdots,
\theta_m \in \R $
\begin{eqnarray} \label{Jiiconv}
J_n&:=&\log \E \exp\Big\{ \i  \frac{1}{n^{1 - \frac{\beta}{\alpha}}} \sum_{j=1}^m \theta_j S_n(\tau_j)  \Big\} \ \to \  J,
\qquad \text{where} \\
J&:=&-\psi(1)\int_{\R_+ \times \R} \Big|\sum_{j=1}^m \theta_j (f(w, \tau_j - u) - f(w, -u)) \Big|^\alpha
\omega \Big(\sum_{j=1}^m \theta_j (f(w, \tau_j - u) - f(w, -u)); \alpha, c^+,c^-\Big)
\frac{\d w \d u}{w^{\alpha -\beta}}. \nn
\end{eqnarray}
We have $J = \log \E \e^{\i  \sum_{j=1}^m \theta_j \Lambda_{\alpha,\beta}(\tau_j)}$
by definition (\ref{ZZ1}) of $\Lambda_{\alpha,\beta}$.
We shall restrict the proof of (\ref{Jiiconv}) to $m= \tau_1 = 1 $, since the general case follows analogously.
Let $V(\theta) $ be defined as in  (\ref{Wchf1}), where $\sigma = 0$.
Then,
\begin{eqnarray*}
J_n
&=&\sum_{s\in \Z} \int_0^\epsilon V\Big(\theta \, \frac{1}{n^{1 - \frac{\beta}{\alpha}}}  \sum_{t=1}^{n} (1-z)^{t-s} \1(t\ge s)\Big) z^\beta \psi(1-z) \d z \\
&=&\sum_{s\le 0} \int_0^\epsilon V(...) z^\beta \psi(1-z) \d z +
\sum_{s=1}^n \int_0^\epsilon V(...) z^\beta \psi(1-z) \d z \\
% + \sum_{s\in \Z} \int_\epsilon^1 V(...) z^\beta \psi(1-z) \d z  \\
& \ =: \ &
J_{n1} + J_{n2}. % + J_{n3},
\end{eqnarray*}
Similarly, split $J = J_1 + J_2$, where
\begin{eqnarray*}
J_1&:=&-|\theta|^\alpha \psi(1) \omega(\theta; \alpha, c^+, c^-)
\int_{-\infty}^0 \d u \int_0^\infty  (f(w,1-u) - f(w,-u))^\alpha
w^{\beta -\alpha} \d w, \nn \\
J_2&:=&-|\theta|^\alpha \psi(1) \omega(\theta; \alpha, c^+, c^-)
\int_0^1 \d u \int_0^\infty  (f(w,u))^\alpha
w^{\beta -\alpha} \d w.
\end{eqnarray*}
To prove \eqref{Jiiconv} we need to show $J_{n1}\to J_1$, $J_{n2} \to J_2$. % $J_{n3} \to 0$.
We shall use the following facts:
\begin{eqnarray} \label{Vlim2}
\lim_{\lambda \to +0} \lambda V\big(\lambda^{-1/\alpha}\theta\big)&=&-|\theta|^\alpha \omega (\theta; \alpha, c^+, c^-), \qquad
\forall \ \theta \in \R
\end{eqnarray}
and
\begin{equation} \label{Vdom}
|V(\theta)| \ \le \  C |\theta|^\alpha, \qquad \forall \, \theta \in \R \qquad (\exists \, C < \infty).
\end{equation}
Here,   (\ref{Vdom}) follows from (\ref{limPi2}),  $\int_{\R} x^2 \pi (\d x) < \infty $ and integration by parts.
To show (\ref{Vlim2}),
let $\chi(x), x \in \R $ be a bounded continuously differentiable function with compact support
and such that $\chi(x) \equiv 1,  |x| \le 1$. Then the l.h.s. of (\ref{Vlim2}) can be rewritten as
\begin{eqnarray*}
\lambda V\big(\lambda^{-1/\alpha}\theta\big)&=&
\int_{\R} (\e^{ {\i} \theta y} - 1 - {\i} \theta y \chi(y)) \pi_\lambda (\d y)
+  \i \theta \mu_{\chi,\lambda},
\end{eqnarray*}
where
$
\pi_\lambda (\d y) := \lambda \pi( \d \lambda^{1/\alpha} y), \
\mu_{\chi,\lambda} :=  \int_{\R} y (\chi(y) -1)  \pi_\lambda ( \d y).
$
The r.h.s. of (\ref{Vlim2}) can be rewritten as
\begin{eqnarray*}
&-|\theta|^\alpha \omega(\theta; \alpha, c^+, c^-) \ = \  V_0(\theta)
\ :=\ \int_{\R} (\e^{ {\i} \theta y} - 1 - {\i} \theta y \chi(y)) \pi_0 (\d y)  + \i \theta \mu_{\chi,0},
\end{eqnarray*}
where
$
\pi_0 (\d y) := -c^+ \d  y^{-\alpha} \1(y>0) + c^- \d (- y)^{-\alpha} \1(y<0), \
\mu_{\chi,0}\ :=\  \int_{\R} y (\chi(y) -1)  \pi_0( \d  y).
$
Let $C_\natural $ be the class of all bounded continuous functions on $\R$ vanishing in a neighborhood of $0$.
According to Sato (1999, Thm. 8.7), relation (\ref{Vlim2}) follows from
\begin{eqnarray}\label{Pi1}
&&\lim_{\lambda \to 0} \int_{\R} f(y) \pi_\lambda (\d y) \ = \  \int_{\R} f(y) \pi_0 (\d y), \quad \forall \, f \in C_\natural, \\
&&
\lim_{\lambda \to 0} \mu_{\chi, \lambda} = \mu_{\chi,0}, \qquad
\lim_{\epsilon \downarrow 0}\lim_{\lambda \to 0} \int_{|y| \le \epsilon} y^2 \pi_\lambda (\d y) \ = \ 0.\label{Pi2}
\end{eqnarray}
Relations (\ref{Pi1}) is immediate from (\ref{limPi2}) while  (\ref{Pi2}) follows from
(\ref{limPi2}) by integration by parts.

\vskip.1cm

Coming back to the proof of (\ref{Jiiconv}),
consider the convergence $J_{n2} \to J_2$.
By change of variables: $n z = w,  n - s +1 = n u $,
$J_{n2}$ can be rewritten as
\begin{eqnarray*}
J_{n2}
&=&
\int_{1/n}^{1} \d u \int_0^{\epsilon n} n^{-\beta} V\Big(\theta n^{\frac{\beta}{\alpha}} \frac{1-  (1- \frac{w}{n})^{[un]}} {w} \Big)
w^\beta \psi\Big(1-\frac{w}{n}\Big) \d w \\
&=&-|\theta|^\alpha \omega(\theta; \alpha, c^+, c^-) \int_0^1 \d u \int_0^{\infty} \Big( \frac{1-  \e^{-wu}} {w} \Big)^\alpha
\kappa_{n2}(\theta; u,w)
w^\beta \psi\Big(1-\frac{w}{n}\Big) \d w,
\end{eqnarray*}
where $\kappa_{n2}(u,w)$ is written as
\begin{eqnarray}
\kappa_{n2}(\theta; u,w)
&:=&-\Big( \frac{1-  \e^{-wu}} {w} \Big)^{-\alpha}  n^{-\beta} \frac{V\Big(\theta n^{\frac{\beta}{\alpha}}w^{-1}(1-  (1- \frac{w}{n})^{[un]}) \Big)}
{|\theta|^\alpha \omega(\theta; \alpha, c^+, c^-)}
\1(n^{-1} < u \le 1, 0< w < \epsilon n) \nn \\
&=&\frac{\lambda_n(u,w) V( \lambda^{-1/\alpha} \theta) }{-|\theta|^\alpha \omega(\theta; \alpha, c^+, c^-)}
 \Big(\frac{1-  (1- \frac{w}{n})^{[un]}} { 1-  \e^{-wu} } \Big)^\alpha \1(n^{-1} < u \le 1, 0< w < \epsilon n) \label{chidef}
\end{eqnarray}
with
\begin{eqnarray*}
\lambda_n(u,w)
&:=&n^{-\beta} \Big(\frac{w}{1-  (1- \frac{w}{n})^{[un]}}\Big)^{\alpha} \  \to \  0
\end{eqnarray*}
for each $u \in (0,1], w >0$ fixed.  Hence and with (\ref{Vlim2}) in mind, it follows that
$\kappa_{n2}(\theta; u,w) \to 1 $ for each $\theta \in \R, u \in (0,1], w >0$ and therefore the convergence
$J_{n2} \to J_2 $ by the dominated convergence theorem provided we establish a dominating bound
\begin{equation} \label{chidom}
|\kappa_{n2}(\theta; u,w)| \ \le \  C
\end{equation}
with $C$ independent of $n, u \in (0,1], w \in (0, \epsilon n).$
From (\ref{Vdom}) it follows that the first
ratio on the r.h.s. of  (\ref{chidef}) is bounded by an absolute constant. Next, for any $0\le x \le 1/2, \, s >0 $ we have
$1- x \ge \e^{-2x}  \, \Longrightarrow \,  (1- x)^s \ge \e^{-2xs}  \, \Longrightarrow \,
1- (1-x)^s \le  2(1- \e^{-xs}) $ and hence $\frac{1-  (1- \frac{w}{n})^{[un]}} { 1-  \e^{-wu} }
\le \frac{1-  (1- \frac{w}{n})^{un}} { 1-  \e^{-wu} } \le 2 $ for any $0\le w \le n/2, \, u>0$ so that the second
ratio on the r.h.s. of  (\ref{chidef}) is also bounded by 2, provided $\epsilon \le 1/2 $. This proves
(\ref{chidom}) and concludes the proof of $J_{n2} \to J_2.$ The proof of the convergence
$J_{n1} \to J_1$ is similar and is omitted.
%Using inequality \eqref{Vdom} it is not difficult to prove that $|J_{n3}| < C n^{\beta-(\alpha-1)}$. Since $\beta-(\alpha-1)<0$, $J_{n3} \to 0$.
This concludes the proof of (\ref{Jiiconv}), or finite-dimensional convergence in (\ref{iiconv}).

To prove the tightness part of (\ref{iiconv}), it suffices to verify the well-known criterion
in Billingsley (1968,  Thm.12.3):  there exists $C>0 $  such that,
for any $n \ge 1 $ and $0 \le \tau < \tau+h \le 1 $
\begin{equation}\label{tightii}
\sup_{u >0} u^\alpha \P\big(n^{\frac{\beta}{\alpha} -1} | S_n(\tau+h) - S_n(\tau) | > u \big) \  < \  C h^{\alpha - \beta},
\end{equation}
where $\alpha - \beta > 1$.
By stationarity of increments of $\{\mathfrak{X}(t)\}$ it suffices to prove (\ref{tightii}) for $\tau = 0, h=1$, in which case it
becomes
\begin{equation}\label{tightii2}
\sup_{u >0} u^\alpha \P\big(| S_n| > u \big) \  < \  C n^{\alpha - \beta}, \qquad S_n := S_n(1).
\end{equation}
The proof of (\ref{tightii2}), below, requires  inequality in (\ref{Mbdd})
for tail probabilities of
stochastic integrals w.r.t. ID random measure.
Let $L^\alpha (\Z \times (-1,1)) $ be the class of measurable functions $g: \Z \times (-1,1) \to \R $  with
$\| g\|^\alpha_\alpha  :=\sum_{s\in \Z} \E |g(s,a)|^\alpha  < \infty.$ Also, introduce
the weak space $L^\alpha_w (\Z \times (-1,1)) $ of measurable functions $g: \Z \times (-1,1) \to \R $  with
$\| g\|^\alpha_{\alpha,w}  :=\sup_{t>0} t^\alpha \sum_{s\in \Z} \P (|g(s,a)|> t)  < \infty.$  Note
$L^\alpha (\Z \times (-1,1)) \subset L^\alpha_w (\Z \times (-1,1)) $  and $\| g\|^\alpha_{\alpha,w} \le  \| g\|^\alpha_{\alpha} $.
Let $\{M_s, s \in \Z\} $ be the random
measure in (\ref{mix}), $M \sim W \sim ID_2(0,\pi)$    with zero mean and the L\'evy measure $\pi $ satisfying
the assumptions in (ii).
It is well-known (see, e.g., Surgailis (1981)) that
the stochastic integral $M(g) :=  \sum_{s \in \Z} \int_{(-1,1)} g(s,a) M_s({\d}a)$ is well-defined for any
$g\in L^p (\Z \times (-1,1)), \, p=1,2$
and satisfies $ \E M^2(g) = C_2\|g \|^2_2, \, \E |M(g)| \le  C_1\|g \|_1 $ for some constants $C_1, C_2 >0$.
The above facts together with Hunt's interpolation theorem, see Reed and Simon (1975, Theorem IX.19) imply
that $M(g) $
extends to all $g \in L^\alpha_w (\Z \times (-1,1)), \, 1 < \alpha < 2 $ and satisfies
the bound
\begin{eqnarray}\label{Mbdd}
\sup_{u>0} u^\alpha \P(|M(g)|> u)&\le&C\|g\|^\alpha_{\alpha,w} \ \le \  C \|g\|^\alpha_{\alpha},
\end{eqnarray}
with some constant $C>0$ depending on $\alpha, C_1, C_2$ only.  Using (\ref{Mbdd}) and the representation $S_n = M(g)$
with $g(s,a) = \sum_{t=1}^n a^{t-s}\1(t\ge s)$ we obtain
\begin{eqnarray*}\label{Sbdd}
\sup_{u >0} u^\alpha \P\big(| S_n| > u \big)&\le&C\sum_{s\le n} \E \Big|\sum_{t=1\vee s}^n a^{t-s} \Big|^\alpha
\ = \ O(n^{\alpha - \beta}),
\end{eqnarray*}
where the last relation easily follows from condition (\ref{bcond}),
see also  Puplinskait\.e and Surgailis (2010, proof of Theorem 3.1). This proves
(\ref{tightii2}) and part (ii).

\vskip.2cm

\noi Proof of (iii).
It suffices to prove that for any $0< \tau_1 < \cdots < \tau_m \le 1, \, \theta_1\in \R, \cdots,
\theta_m \in \R $
\begin{eqnarray} \label{Jiiiconv}
J_n\ :=\ \log \E \exp\Big\{ \i  \frac{1}{n^{1/(1+\beta)}} \sum_{j=1}^m \theta_j S_n(\tau_j)  \Big\}
&\to&
J \ := \ \log \E \exp\{ \i \sum_{j=1}^m \theta_j L_{1+ \beta}(\tau_j)  \}.
\end{eqnarray}
Similarly as in (i)-(ii), we shall restrict the proof of (\ref{Jiiiconv}) to the case
$m= 1$ since the general case follows analogously.
Then
\begin{eqnarray*}
J_n
&=&\sum_{s\in \Z} \int_0^\epsilon V\Big(n^{-1/(1+\beta)} \theta \sum_{t=1}^{[n\tau]} (1-z)^{t-s} \1(t\ge s)\Big) z^\beta \psi(1-z) \d z \ = \
J_{n1} + J_{n2},
\end{eqnarray*}
where $J_{n1} := \sum_{s\le 0} \int_0^\epsilon V(\cdots) z^\beta \psi(1-z) \d z, \, J_{n2} := \sum_{s=1}^{[n\tau]}
\int_0^\epsilon V(\cdots) z^\beta \psi(1-z) \d z. $
Let $\theta > 0$. By the change of variables: $n^{1/(1+\beta)}  z = \theta/y,  [n\tau] - s +1 = n u $,
$J_{n2}$ can be rewritten as
\begin{eqnarray}
J_{n2}
&=&\sum_{s=1}^{[n\tau]} \int_0^\epsilon V\Big(\frac{\theta(1 - (1-z)^{[n\tau] - s +1})}{n^{1/(1+\beta)}  z}
\Big) z^\beta \psi(1-z) \d z \nn \\
&=&\theta^{1+ \beta} \int_0^\tau \d u \int_0^\infty \frac{ \d y}{y^{\beta +2}}
 V\Big(y(1-  (1- \frac{\theta}{n^{1/(1+\beta)} y})^{[un]}) \Big) \psi \Big(1-\frac{\theta}{n^{1/(1+\beta)} y}\Big ) \1_n(\theta; y,u), \label{Jn2iii}
\end{eqnarray}
where $\1_n(\theta; y,u) := \1( 1/n < u < [n\tau]/n], y > \theta \epsilon^{-1} n^{-1/(1+ \beta)}) \to \1 (0< u < \tau, y >0).$  As
$(1- \frac{\theta}{n^{1/(1+\beta)}y})^{un} \to 0$ for any $u, y >0$  due to $n /n^{1/(1+ \beta)} \to \infty, $ we see that
the integrand in (\ref{Jn2iii}) tends to $y^{-\beta -2} V(y)\psi(1). $  We will soon prove  that this passage to the limit
under the sign of the integral in  (\ref{Jn2iii}) is legitimate. Therefore,
\begin{eqnarray}
J_{n2}
&\to&J \ := \  \tau |\theta|^{1+ \beta}\psi(1) \int_0^\infty V(y) y^{-\beta -2} \d y  \ = \  -\tau |\theta|^{1+ \beta} \psi(1) \omega(\theta; 1+ \beta,
\pi_\beta^{-}, \pi^{+}_\beta), \label{Jiii} \\
\pi^+_\beta&:=&\frac{1}{1+\beta}\int_{0}^{\infty} x^{1+ \beta} \pi (\d x), \qquad \pi^-_\beta \ :=\ \frac{1}{1+\beta}\int_{-\infty}^{0} |x|^{1+ \beta} \pi (\d x), \nn
\end{eqnarray}
and the last equality in (\ref{Jiii}) follows from  the definition of $V(y)$ and
Ibragimov and Linnik (1971, Thm. 2.2.2).

For justification of the above passage to the limit, note that the function
$V(y) = \int_{\R} (\e^{\i y x} - 1 - \i yx)\pi (\d x)$ satisfies
$|V(y)|\le V_1(y) + V_2(y),$ where
$V_1(y) :=  y^2 \int_{|x| \le 1/|y|} x^2 \pi (\d x), \ V_2(y) :=  2|y|\int_{|x| > 1/|y|} |x| \pi (\d x). $
We have
\begin{eqnarray*}
\int_0^\infty (V_1(y)+ V_2(y))y^{-\beta -2}\d y&\le&\int_\R x^2 \pi (\d x) \int_0^{1/|x|} y^{-\beta} \d y
+ 2\int_\R |x| \pi (\d x) \int_{1/|x|}^\infty y^{-1-\beta} \d y \\
&\le&C \int_{\R} |x|^{1+ \beta} \pi (\d x) \ < \ \infty.
\end{eqnarray*}
Next, $\sup_{1/2 \le c \le 1} V_1(cy) \le y^2  \int_{|x| \le 2/|y|} x^2 \pi (\d x) =: \bar V_1(y), \,
\sup_{1/2 \le c \le 1} V_2(cy) \le V_2(y) $ and $\int_0^\infty \bar V_1(y)y^{-\beta -2}\d y < \infty $.
Denote $\zeta_n(\theta; y,u) := (
1- \frac{\theta}{n^{1/(1+\beta)}y})^{[un]} $. Then $ \zeta_n(\theta; y,u) \ge 0$ and we split the integral
in (\ref{Jn2iii}) into two parts corresponding to $\zeta_n(\theta; y,u) \le 1/2 $ and
$\zeta_n(\theta; y,u) > 1/2 $, viz.,  $J_{n2} = J_{n2}^+ + J_{n2}^-$, where
\begin{eqnarray*}
J^+_{n2}
&:=&\theta^{1+ \beta} \int_0^\tau \d u \int_0^\infty y^{-\beta -2} \d y
V\Big(y(1-  \zeta_n(\theta; y,u)) \Big) \psi \Big(1-\frac{\theta}{n^{1/(1+\beta)} y}\Big ) \1 (\zeta_n(\theta; y,u) \le 1/2)   \1_n(\theta, y,u), \\
J^-_{n2}
&:=&\theta^{1+ \beta} \int_0^\tau \d u \int_0^\infty y^{-\beta -2} \d y
V\Big(y(1-  \zeta_n(\theta; y,u)) \Big) \psi \Big(1-\frac{\theta}{n^{1/(1+\beta)} y}\Big ) \1 (\zeta_n(\theta; y,u) > 1/2)   \1_n(\theta; y, u).
\end{eqnarray*}
Since $\big|V\big(y(1-  \zeta_n(\theta; y,u)) \big)  \1 (\zeta_n(\theta; y,u)  \le 1/2)\big| \le \bar V_1(y) + V_2(y) $ is bounded
by integrable function (see above), so $J^+_{n2} \to J$ by the dominated convergence theorem. It remains
to prove $J_{n2}^- \to 0.$ From inequalities $ 1- x \le \e^{-x}\ (x >0) $ and $[un] \ge un/2 \ (u > 1/n)$
it follows that $\zeta_n(\theta; y,u) \le \e^{ - \theta un/ 2n^{1/(1+\beta)}y} $ and hence
$\1 (\zeta_n(\theta; y,u) > 1/2) \le \1(  \e^{ - \theta un/ 2n^{1/(1+\beta)}y} > 1/2 )  =
\1 ( (u/y) < c_1 n^{-\gamma}), $  where $\gamma := \beta/(1+ \beta) >0, \, c_1 := (2 \log 2)/\theta $. Without loss
of generality, we can assume that $1 < \alpha < 1+ \beta$ in  (\ref{pimom}).
Condition (\ref{pimom}) implies
\begin{eqnarray*}
|V(y)|
&\le&\int_{|xy| \le 1} |y x|^{\alpha} \pi (\d x) + 2 \int_{|y x| >1} |y x|^{\alpha} \pi (\d x)  \ \le \
C|y|^{\alpha}, \qquad \forall \, y \in \R.
\end{eqnarray*}
Hence
\begin{eqnarray*}
|J^-_{n2}|
&\le&C\int_0^\tau \d u \int_0^\infty  \1 \big( \frac u y < c_1 n^{-\gamma}\big) \frac{ \d y}{y^{2+ \beta - \alpha}} \ \le \ K
n^{-\gamma(1+ \beta-\alpha)}  \   \to \  0,
\end{eqnarray*}
where $K := C \int_0^\tau u^{\alpha -1-\beta} \d u < \infty $.
This proves
$J_{n2} \to J$, or
(\ref{Jiii}). The proof of $J_{n1} \to 0$ follows similarly and hence is omitted.

\vskip.2cm

\noindent Proof of (iv). The proof of finite-dimensional convergence is similar
to Puplinskait\.e and Surgailis (2010, proof of Thm. 3.1 (ii)). Below, we present the
proof of the one-dimensional convergence of $n^{-1/2} S_n = n^{-1/2} \sum_{t=1}^n \mathfrak{X}(t) $ towards  ${\cal N}(0, \sigma^2_\Phi)$ with $\sigma^2_\Phi >0$ given in
(\ref{sigmaPhi}) below.  The convergence of general finite-dimensional distributions follows analogously.
Similarly as above, consider $J_n := \log \E \exp\{\i \theta n^{-1/2} S_n\} =  J_{n1} + J_{n2}$, where
$J_{n1} := \sum_{s \le 0} \E V \big(\theta n^{-1/2} \sum_{t=1}^n a^{t-s} \big), \
J_{n2}:=\sum_{s=1}^n \E V \big(\theta n^{-1/2} \sum_{t=s}^n a^{t-s} \big). $  Let $\tilde \Phi(\d z) := \Phi(\d (1-z)), z\in (0,2)$.
We have
\begin{eqnarray*}
J_{n2}&=&\sum_{k=1}^n \int_{(0,2)} V \big(\theta \frac{ 1- (1-z)^k}{z n^{1/2}} \big)
%\phi(1-z)
\tilde \Phi(\d z) \\
&=&-
\theta^2 \sigma^2_W \, n^{-1} \sum_{k=1}^n \int_{(0,2)} (1 - (1-z)^k )^2 z^{-2}
 \kappa_n(\theta;k,z) \tilde \Phi(\d z),
%\phi(1-z) \d z,
\end{eqnarray*}
where  $\kappa_n(\theta; k,z):=  \kappa \big(\theta \frac{ 1- (1-z)^k}{z n^{1/2}} \big)$ and the function
$\kappa(y):=-\frac{V(y)}{\sigma^2_W y^2}$ satisfies $\lim_{y \to 0} \kappa(y) = 1, \, \sup_{y \in \R} |\kappa(y)| < \infty$.
These facts together with $\beta > 1$ imply $ n^{-1} \sum_{k=1}^n \int_{(0,2)} (1 - (1-z)^k )^2 z^{-2} \kappa_n(\theta;k,z)
%\phi(1-z) \d z
\Phi(\d z)
\to \int_{(0,2)} z^{-2} \tilde \Phi(\d z)$
%\phi(1-z) \d z $
and hence $ J_{n2} \to - (1/ 2) \theta^2 \sigma^2_\Phi$, with
\begin{equation}\label{sigmaPhi}
\sigma_\Phi^2 := 2\sigma^2_W \int_{(0,2)} z^{-2} \tilde \Phi(\d z)% \phi(1-z) \d z
\ = \ 2\sigma^2_W \E (1-a)^{-2}.
\end{equation}
The proof of $J_{n1} \to 0$ follows similarly (see Puplinskait\.e and Surgailis (2010) for details). This proves
(\ref{ivconv}).

Let us prove the tightness part in (iv). It suffices to show the bound
\begin{equation}\label{S4bdd}
\E S^4_n \ \le \ C n^{2}.
\end{equation}
We have $S_n = M(g)$, where
$M $ is the stochastic integral discussed in the proof of (ii) above and $g \equiv g(s,a) = \sum_{t=1}^n a^{t-s} \1 (t\ge s) \in
L^2(\Z \times (-1,1)).$
Then
$\E M^4 (g) = {\rm cum}_4 (M(g)) + 3 (\E M^2(g) )^2, $ where  $ \E M^2(g) = \E S^2_n $ satisfies $ \E S^2_n \le C n $
(the last fact follows by a similar argument as above).
Hence,
$(\E M^2(g) )^2  \le C n^2 $ in agreement with (\ref{S4bdd}). It remains to evaluate the 4th cumulant
${\rm cum}_4 (S_n) =
{\rm cum}_4 (M(g))=    \pi_4 \sum_{s\in \Z} \E g^4(s,a)$, where $\pi_4 := \int_{\R} x^4 \pi (\d x). $
Then  ${\rm cum}_4 (S_n)=\pi_4 (L_{n1} + L_{n2}), $ where
\begin{eqnarray*}
L_{n1}&:=&\sum_{s\le 0} \E \Big(\sum_{t=1}^n a^{t-s}\Big)^4, \qquad L_{n2}\ :=\ \sum_{s=1}^n \E \Big(\sum_{t=s}^n a^{t-s}\Big)^4.
\end{eqnarray*}
We have
\begin{eqnarray*}
L_{n2}&\le&n\sum_{k=1}^n \E \Big|\sum_{t=0}^k a^t\Big|^3 \ \le \  n \sum_{k=1}^n \E \big[\frac{1}{|1-a|^3} \big]
%\int_0^1  z^{ \beta-3} \psi(1-z) \d z
\ \le \  Cn^2
\end{eqnarray*}
since $\beta > 2 $. Similarly,
\begin{eqnarray*}
L_{n1}&\le&n^2\sum_{s\le 0} \E \Big(\sum_{t=1}^n a^{t-s}\Big)^2 \ \le \  n^2 \E \big[ \frac{1}{(1-a^2)(1-a)^2} \big]
%\int_0^1 \frac{z^{\beta -2} \psi(1-z) \d z}{1 - (1-z)^{ 2 }}
\ \le \  Cn^2.
\end{eqnarray*}
This proves (\ref{S4bdd}) and part (iv). Theorem \ref{sums} is proved. \hfill $\Box$

\medskip

%%%%%%%%%%%%%%%%%%%%%%%%%%%%%%%%%%%%%%%%%%%%%%%%%%%%%%%
\section{Disaggregation}
%%%%%%%%%%%%%%%%%%%%%%%%%%%%%%%%%%%%%%%%%%%%%%%%%%%%%%%

%
Following Leipus et al. (2006), let us define an estimator of $\phi$, the density of the mixing distribution $\Phi$. Differently
from the last paper,
we shall assume below that the variance $\sigma^{-2}_W$ is not
necessary known.
%$\Phi $ is concentrated on the interval $(0,1)$. % and has a density $\phi(x), 0< x < 1$.
%The last restriction is not very important but leads to some differences in the definition of the estimator.
Its
starting point is the equality (\ref{covX}), implying
\begin{equation} \label{covX2}
\sigma^{-2}_W (r(k) -  r(k+2)) \ = \ \int_{-1}^1  x^k  \phi (x) \d x, \qquad k = 0,1, \cdots,
\end{equation}
where $r(k) = {\rm Cov}(\mathfrak{X}(k),  \mathfrak{X}(0))$ and $\sigma^2_W=\var(W) = r(0) - r(2)$.
The l.h.s. of (\ref{covX2}), hence  the integrals on the r.h.s. of (\ref{covX2}),
or moments of $\Phi$,
can be estimated from the observed
sample, leading to the problem of recovering the density from its moments, as explained below.

For a given $q>-1$,  consider a finite measure on $(-1,1)$ having
density $w^{(q)}(x):=(1-x^2)^{q}$. Let $L_2(w^{(q)})$
be the space of functions $h: (-1,1) \to \R$ which are square integrable with respect
to  this measure. Denote by $\set{G^{(q)}_n,\; n=0,1,\cdots}$
the orthonormal basis in $L_2(w^{(q)})$ consisting of normalized
Gegenbauer polynomials  $G^{(q)}_n(x) = \sum_{j=0}^n g^{(q)}_{n,j} x^j$ 
%\begin{equation}
%\label{jacobi}
%  G^{(q)}_n(x)= \frac{
% {g_n}^{-1/2}}{\Gamma(q+1/2) }
% \sum_{m=0}^{|n/2|}(-1)^m   \frac{\Gamma(q+1/2+n-m)}{\Gamma(m+1)\Gamma(n-2m+1)}\
% (2x)^{n-2m}.
%\end{equation}
with coefficients 
\begin{equation}
\label{jcoeff}
g^{(q)}_{n,n-2m} = (-1)^m \frac{
(g_n)^{-1/2}}{\Gamma(q+1/2) }
  \frac{2^{n-2m}\Gamma(q+1/2+n-m)}{\Gamma(m+1)\Gamma(n-2m+1) } \qquad
\mathrm{for} \quad 0\leq m \leq \left[n/2\right],
\end{equation}
where $g_n := \frac{\pi}{2^{2q}} \frac{\Gamma (n + 2q +1)}{\Gamma^2 (q+ 1/2) \Gamma (n+q + 1/2)}$, 
see Abramovitz and Stegun (1965, (22.3.4)), also Leipus et al. (2006, (B.4)).  
Thus,
\begin{equation} \label{orthoJ}
  \int_{-1}^1 G^{(q)}_j(x) G^{(q)}_k(x) w^{(q)}(x) \d x \
  = \  \left\{ \begin{array}{ll}
1&\quad \hbox{if}\;j=k,\\
0&\quad \hbox{if}\;j\neq k.
\end{array}\right.
 \end{equation}
Any function $h \in L_2(w^{(q)})$ can be expanded in Gegenbauer polynomials:
\begin{equation}\label{hexp}
h(x)=\sum_{k=0}^\infty h_k G^{(q)}_k(x) \quad \mathrm{with}
\quad h_k=\int_{-1}^1 h(x)G^{(q)}_k(x) w^{(q)}(x) \d
x =\sum_{j=0}^k g^{(q)}_{k, j}  \int_{-1}^1 h(x)x^{j} w^{(q)}(x) \d
x.
\end{equation}
Below, we call (\ref{hexp}) the $q$-Gegenbauer expansion of $h$.

Consider the function
\begin{equation}\label{zphi}
\zeta(x) := \frac{\phi(x)}{(1-x^2)^{q} }, \qquad \text{with} \quad \int_{-1}^1  \zeta (x) (1-x^2)^{q} \d x = \int_{-1}^1  \phi (x)  \d x = 1.
\end{equation}
Under the condition
\begin{equation}
\label{phicond} \int_{-1}^1
\frac{\phi(x)^2}{(1-x^2)^{q
}} \d x <\infty,
\end{equation}
the function $\zeta$ in \eqref{zphi} belongs to $L_2(w^{(q)})$, and  has a $q-$Gegenbauer expansion with coefficients
\begin{equation}
 \label{zetacoef}
 \zeta_k=
 \sum_{j=0}^k g^{(q)}_{k,j} \int_{-1}^1 \phi(x)
  x^j \d x = \frac{1}{\sigma_W^2}
\sum_{j=0}^k g^{(q)}_{k,j} \left(r{(j)}-r(j+2)\right),
  \qquad k= 0,1, \cdots;
\end{equation}
see (\ref{covX2}). Equations (\ref{hexp}), (\ref{zetacoef})
lead to the following estimates of the function $\zeta(x)$:
\begin{eqnarray} \label{zetaest}
\widehat \zeta_n(x)&:=&\sum_{k=0}^{K_n} \widehat \zeta_{n,k} G^{(q)}_k(x),
\qquad
\widetilde \zeta_n(x)\ :=\  \sum_{k=0}^{K_n} \widetilde \zeta_{n,k} G^{(q)}_k(x),
\end{eqnarray}
where $K_n, n \in \N^* $ is a nondecreasing sequence tending to infinity at a rate which is discussed below,
and
\begin{eqnarray}\label{zetacoef_m_est}
\widehat \zeta_{n,k}&:=&\frac{1}{\widehat \sigma_W^2}
\sum_{j=0}^k g^{(q)}_{k,j}
(\widehat r_n(j)-\widehat r_n(j+2)),
\qquad
\widetilde \zeta_{n,k}\ :=\  \frac{1}{\sigma_W^2}
\sum_{j=0}^k g^{(q)}_{k,j}
(\widehat r_n(j)-\widehat r_n(j+2))
\end{eqnarray}
are natural estimates of the $\zeta_k$'s in (\ref{zetacoef}) in the case when
$\sigma^2_W$ is unknown or known, respectively. Here and below,
\begin{equation} \label{covest}
\overline{\mathfrak{X}} := \frac 1 n \sum_{k=1}^n \mathfrak{X}(k), \qquad \widehat r_n(j):= \frac{1}{n} \sum_{i=1}^{n-j}
\big(\mathfrak{X}(i) - \overline{\mathfrak{X}}\big) \big(\mathfrak{X}(i+j)- \overline{\mathfrak{X}}\big), \quad j=0,1, \cdots, n
\end{equation}
are the sample mean and the sample covariance, respectively, and the estimate
of $\sigma^2_W = r(0) - r(2)$
is defined as $$\widehat \sigma_W^2 :=  \widehat r_n(0) -  \widehat r_n(2).  $$
The corresponding estimators of $\phi (x)$ is constructed following relation (\ref{zphi}):
\begin{eqnarray} \label{phiestnew}
\widehat \phi_n(x)&:=&\widehat \zeta_n(x)(1-x^2)^{q}, \qquad \widetilde \phi_n(x)\ :=\  \widetilde \zeta_n(x)(1-x^2)^{q}.
\end{eqnarray}
The above estimators
were essentially constructed in Leipus et al. (2006) and
Celov et al. (2010). The modifications in  (\ref{phiestnew}) differ from the original ones in the above mentioned
papers by the choice of a more natural estimate (\ref{covest}) of the covariance function
$r(j)$, which allows for non-centered observations and makes both  estimators
in (\ref{phiestnew})
location and scale invariant. Note also that
the first estimator in (\ref{phiestnew}) satisfies $\int_{-1}^1  \widehat \phi_n(x) \d x = 1, $ while the second
one does not have this property and can be used only if $\sigma^2_W$ is known.

\begin{prop}\label{disagg}
Let $(\mathfrak{X}(t))$ be an aggregated process in (\ref{mix}) with finite 4th moment $\E\mathfrak{X}(0)^4 < \infty $ and
$M \sim W \sim ID(\mu , \sigma, \pi)$.  Assume that the mixing
density $\phi(x) $  satisfies conditions  (\ref{moment})
and (\ref{phicond}), with some $q>-1$.
Let $ \widetilde \zeta_n(x)$ be the estimator of $\zeta(x)$ as defined in \eqref{zetaest}, where
$K_n$ satisfy
\begin{equation}\label{Kn}
K_n= [\gamma \log n]\quad \mathrm{with}\quad 0<\gamma < (2\log
(1+\sqrt{2}))^{-1},
\end{equation}
Then
\begin{equation}
  \label{disconv}
\int_{-1}^1\E{(\widetilde
  \zeta_n(x)-\zeta(x))^2}(1-x^2)^{q} \d x \ \to \  0.
\end{equation}
\end{prop}

\noi {\it Proof.} Denote $v_n$ the l.h.s. of (\ref{disconv}).
From the orthonormality property (\ref{orthoJ}), similarly as in Leipus et al. (2006, (3.3)),
\begin{equation}\label{vn}
v_n  \ = \ \sum_{k=0}^{K_n} \E (\widetilde \zeta_{n,k} - \zeta_k)^2
 + \sum_{k=K_n+1}^\infty \zeta^2_k,
\end{equation}
where the second sum on the r.h.s. tends to 0. By the location invariance mentioned above, w.l.g. we can
assume below that $  \E\mathfrak{X}(t) = 0$.
Let $\widehat r^\circ_n(j):= \frac{1}{n} \sum_{i=1}^{n-j}
\mathfrak{X}(i)\mathfrak{X}(i+j), \ 0\le j <n, $ then $\E \widehat r^\circ_n(j) - r(j) = (j/n) r(j) $ and
\begin{eqnarray}
\E \big\{\widetilde \zeta_{n,k} - \zeta_k\big\}^2
&=&\sigma^{-4}_W \E\Big\{\sum_{j=0}^k g^{(q)}_{k,j} \big( \widehat r_n(j) - \widehat r_n(j+2)-r(j) +r(j+2)\big) \Big\}^2 \nonumber \\
&=&\sigma^{-4}_W \E\Big\{\sum_{j=0}^k g^{(q)}_{k,j} \Big( \widehat r^\circ_n(j) - \widehat r^\circ_n(j+2)-r(j) +r(j+2)+
2 n^{-1}\overline{\mathfrak{X}}^2  \nonumber \\
&&\hskip1cm -\ n^{-1} \overline{\mathfrak{X}}\big[\mathfrak{X}(n-j-1)+
\mathfrak{X}(n-j) +\mathfrak{X}(j+1)+ \mathfrak{X}(j+2)\big]\Big) \Big\}^2 \nonumber \\
&\le&C k \big(\max_{0\le j \le k} |g^{(q)}_{k,j}| \big)^2 \sum_{j=0}^k \Big(\frac{j^2}{n^2} + {\rm Var}(\widehat r^\circ_n(j) - \widehat r^\circ_n(j+2))+\frac{C}{n^2}\Big),\label{zetan}
\end{eqnarray}
where we used the trivial bound $\E \overline{\mathfrak{X}}^4<C$.

The rest of the proof of Proposition \ref{disagg}
follows from (\ref{vn}), (\ref{zetan}), Lemmas \ref{lem:sigmak} below
and the following  bound on the Gegenbauer coefficients
%obtained by Leipus et al. (2006, Lemma 5)
$$
\max_{0\le j \le n }{|g^{(q)}_{n,j}|} \leq  C
n^{11/2} \e^{n\beta} \quad \hbox  {with}\quad  \beta :=
\log(1+\sqrt{2}),
$$
% for
%any $n$.
\noi obtained in Leipus et al. (2006, Lemma 5). See Leipus et al. (2006, pp.2552-2553) for other details. \hfill $\Box$

\smallskip

\noi Lemma \ref{lem:sigmak}
generalizes  (Leipus et al., 2006, Lemma 4) for a non-Gaussian aggregated process
with finite 4th moment.

\begin{lemma}\label{lem:sigmak}
Let $\{\mathfrak{X}(t)\}$ be an aggregated
process in (\ref{mix}) with $\E \mathfrak{X}(0)^4 < \infty, \, \E \mathfrak{X}(0) = 0 $.
There exists a  constant $C>0$ independent  of $n, k$ and such that
\begin{equation}\label{eq:varvar}
 \var(\widehat r^\circ_n(k)-\widehat r^\circ_n(k+2))\leq\frac{C}{n}.
\end{equation}

\end{lemma}

\noi {\it Proof.}
Let $D(k):= \mathfrak{X}(k)-\mathfrak{X}(k+2).$ Similarly as in Leipus et al. (2006, p.2560),
$$
 \var(\widehat r^\circ_n(k)-\widehat r^\circ_n(k+2)) \ \le \  Cn^{-2} \Big( \var \big(\sum_{j=1}^{n-k-2} \mathfrak{X}(j)D(j+k)\big) + 1\Big).
$$
Here, $\var \big(\sum_{j=1}^{n-k-2} \mathfrak{X}(j) D(j+k)\big) = \sum_{j,l=1}^{n-k-2} \cov \big(\mathfrak{X}(j) D(j+k),
\mathfrak{X}(l) D(l+k)\big)$, where
\begin{eqnarray*}
\cov(\mathfrak{X}(j)D(j+k),\mathfrak{X}(l) D(l+k))&=&{\rm Cum}(\mathfrak{X}(j),D(j+k),\mathfrak{X}(l), D(l+k))\\
&+&\E[\mathfrak{X}(j)\mathfrak{X}(l)]\E[D(j+k)D(l+k)] + \E[\mathfrak{X}(j)D(k+l) ]\E[\mathfrak{X}(l)D(j+k)].
\end{eqnarray*}
The two last terms in the above representation of the covariance are estimated in Leipus et al. (2006). Hence
the lemma follows from
\begin{equation}\label{cumest}
\sum_{j,l=1}^{n-k-2}{\rm Cum}(\mathfrak{X}(j),D(j+k),\mathfrak{X}(l), D(l+k))\ \leq\ Cn.
\end{equation}
We have for $k_1,k_2 \ge 0, l \ge j$
\begin{eqnarray*}
{\rm Cum}(\mathfrak{X}(j),\mathfrak{X}(j+k_1),\mathfrak{X}(l), \mathfrak{X}(l+k_2))
&=&\pi_4 \E \Big[ \sum_{s \le j} a^{j-s} a^{j-s+k_1}  a^{l-s} a^{l-s+k_2} \Big] \\
&=&\pi_4 \E \Big[\frac{a^{k_1 + k_2 +2(l-j)}}{1- a^4} \Big]
\end{eqnarray*}
and hence
\begin{eqnarray*}
c_{j,l,k} \ := \ {\rm Cum}(\mathfrak{X}(j),D(j+k),\mathfrak{X}(l), D(l+k))
&=&\pi_4 \E \Big[\frac{a^{2k +2(l-j)} (1-a^2)}{1+ a^2} \Big]
\end{eqnarray*}
where $\pi_4 := \int_{\R} x^4 \pi({\d}x) $. Then
\begin{eqnarray*}
\sum_{j,l=1}^{n-k-2}|c_{j,l,k}|
&\le&C\sum_{1\le j\le l \le n} \E \big[\frac{(1-a^2)}{1+a^2} |a|^{2(l-j)} \big] \\
&\le&C\sum_{1\le j \le n} \E \big[\frac{1}{1+a^2} \big] \ \le \  Cn,
\end{eqnarray*}
proving (\ref{cumest}) and the lemma, too. \hfill $\Box$

 The main result of this sec. is the following theorem.

\begin{theorem}\label{disagg2}
Let $\{\mathfrak{X}(t)\}$, $\phi(x) $ and $K_n$ satisfy the conditions of Proposition \ref{disagg},
and $ \widehat \phi_n(x), \widetilde \phi_n(x)$ be the estimators of $\phi(x)$ as defined in \eqref{phiestnew}.
Then
\begin{equation}
  \label{disconv1}
\int_{-1}^1 \frac{(\widehat
  \phi_n(x)-\phi(x))^2}{(1-x^2)^{q}} \d x \ \to_p \  0 \qquad \text{and} \qquad \int_{-1}^1 \frac{\E (\widetilde
  \phi_n(x)-\phi(x))^2}{(1-x^2)^{q}} \d x \ \to \  0.
\end{equation}
\end{theorem}

\noi {\it Proof.} The second relation in (\ref{disconv1}) is immediate from \eqref{phiestnew} and
(\ref{disconv}). Next,
$$\widehat \phi_n(x)-\phi(x) = \frac{\sigma^2_W}{ \widehat \sigma^2_W} \big(\widetilde \phi_n(x)-\phi(x)\big) + \phi(x) \big( \frac{\sigma^2_W}{ \widehat \sigma^2_W} -1\big),$$ where
$$
 \widehat \sigma^2_W = \widehat r_n(0) - \widehat r_n(2) = (g_{0,0}^{(q)})^{-1}  \sigma^2_W \widetilde \zeta_{n,0} = \sigma^2_W \int_{-1}^1   \widetilde \zeta_n(x) (1-x^2) ^{q} \d x,
$$
see (\ref{orthoJ}), (\ref{hexp}),  (\ref{zetaest}), (\ref{zetacoef_m_est}). Hence the first relation in (\ref{disconv1})
follows from the second one and the fact that $\widehat \sigma^2_W - \sigma^2_W \to_p \  0$. We have
\begin{eqnarray*}
\E(\widehat \sigma^2_W - \sigma^2_W)^2&=&\sigma^{4}_W \E  \Big (\int_{-1}^1 (\widetilde \zeta_n(x)-\zeta(x)) (1-x^2)^{q} \d x \Big )^2\\
&\leq&\sigma^{4}_W \E  \Big (\int_{-1}^1 (\widetilde \zeta_n(x)-\zeta(x))^2 (1-x^2)^{q} \d x \int_{-1}^1  (1-x^2)^{q} \d x \Big )\\
&=& \sigma^4_W 2^{2q+1}\frac{\Gamma(q+1)^2}{\Gamma(2q+1)} \int_{-1}^1 \E(\widetilde \zeta_n(x)-\zeta(x))^2 (1-x^2)^{q} \d x \to \  0, \ \text{ as $n\to \infty,$}
\end{eqnarray*}
see (\ref{disconv}).
Theorem \ref{disagg2} is proved. \hfill $\Box$

\begin{rem}\label{rem4.2}
{\rm An interesting open question is asymptotic normality of the mixture density estimators in (\ref{phiestnew})
for non-Gaussian process  $\{\mathfrak{X}(t)\}$   (\ref{mix}),
extending Theorem 2.1 in Celov et al. (2010). The proof of the last result relies on a central limit
theorem for quadratic forms of moving-average processes due to
Bhansali et al. (2007). Generalizing this theorem to mixed ID moving averages is an open problem at this moment.
}

\end{rem}

%\bigskip

%\begin{example}

\noi {\bf A simulation study.} We illustrate the performance of the estimator $\widehat\phi_n$
in \eqref{phiestnew} from aggregated processes with Gamma and Gaussian innovations.
%simulated with Gamma noise.
Write $\xi \sim {\rm Gamma}(a,b)$ if $\xi$ has gamma distribution with density
proportional to $x^{a-1} \e^{-x/b} \1_{(0,\infty)}(x) $, with  mean $ab$ and  variance $ab^2$.
It is well-known that $\xi \sim  \text{Gamma} (a,b) $ is
ID and  $\E \e^{\i \theta \xi} = (1- \i \theta b)^{-a} = \exp\{ \int_0^\infty (1-\e^{\i \theta x}) \d \Pi^+(x)\},
\Pi^+ (x) := a \int_x^\infty y^{-1} \e^{-y/b} \d y,  x >0.  $ The statistics
$\widehat\phi_n$ is computed for the aggregated process  $\mathfrak{X}_N(t) = \sum_{i=1}^N X_i^{(N)}(t), 1\le t \le n$ with
$N = 5,000$ and $\{X_i^{(N)}(t)\}$ simulated according to the AR(1) equations in (\ref{AR}).
%The innovation distribution $\vep^{(N)} = \xi^{(N)} - \E \xi^{(N)},
%\xi^{(N)} \sim  \text{Gamma} (1/N,1) $ satisfies   (\ref{DID}) with $\to_{\rm d}$ replaced by $=_{\rm d}$ and
%$W =_{\rm d} \xi - \E \xi, \, \xi \sim  \text{Gamma} (1, 1)$, the centered Gamma(1,1) distribution.
We consider two cases of the noise distribution %in $\epsilon_N(t)$
in \eqref{AR}:
\begin{eqnarray}
\varepsilon^{(N)}(t)&\sim&\text{Gamma} (1/N,1) - 1/N, \label{notcl} \\
%    G_N-\frac 1N \quad \text{where} \quad G_N \sim \Gamma(\frac1N,1)
%  \end{equation}
%The aggregated process  is non gaussian (as shown Figure\ref{fig:traj})
\varepsilon^{(N)}(t)&\sim&{\cal  N}(0,1/N).
%has the same distribution  as
%  \begin{equation}
\label{gamtcl}
\end{eqnarray}
In our simulations, we take the mixing
distribution with density %$\phi$ of the form
\begin{equation}
\phi(x) \propto (1+x) (1-x)^\beta \1_{(-1,1)}(x),\label{exdens}
\end{equation}
with $\beta $ taking values $0.25, 0.75$ and $1.25$. Thus, for $\beta = 0.25, 0.75$
the aggregated process has covariance long memory
and for $\beta =1.25$ it has covariance short memory in both cases
\eqref{notcl} and \eqref{gamtcl}. The simulated trajectory with Gamma innovations \eqref{notcl}
shown in Figure 1
clearly indicates that this process is nongaussian. The L\'evy measure of  \eqref{notcl}
satisfies the asymptotics in \eqref{limPi2} with $\alpha =0$ up to a logarithmic factor.
Following the proof of Theorem \ref{sums} (iii), it can be easily shown that partial sums of the limit aggregated process in the case
\eqref{notcl} tends to a $(1+\beta)-$stable L\'evy process for any $0<\beta <1$, thus also
for $\beta = 0.25 $ and $0.75$.

The estimate $\widehat\phi_n$ strongly depends on $q$ and $K_n$. For $\phi $
in (\ref{exdens}), condition (\ref{phicond}) is satisfied with any
$-1< q < 1+ 2\beta$. In particular, $q<1$ ensures this condition for arbitrary
$\beta>0$, which is generally unknown. 

Figure \ref{fig:real} illustrates the behavior of the estimate
$\widehat\phi_n$  when the distribution of the noise is given  by
\eqref{notcl}.  Here, the parameter $q=0.5$ is fixed. This figure clearly shows the presence of a strong bias for
smaller values of $K_n = 0,1,2$  and an increase in the variance for $K_n=3,4$. 
Figure 1 also suggests that the accuracy of the estimate decreases 
with $\beta$, or with the memory increasing in the aggregated process.  

Figures  \ref{fig:miseGam} and \ref{fig:misenotcl} represent integrated MISE of $\widehat\phi_n$ estimated by a Monte Carlo procedure
with 500 replications, for models 
\eqref{notcl} - \eqref{exdens} and different values of parameters $q$ and $\beta $. 
While the optimal choice of $q$ (minimizing the integrated MISE in
(\ref{disconv1})) is not clear, Figures \ref{fig:miseGam} and \ref{fig:misenotcl} suggest that the ``optimal'' choice of $q$
might be close to (unknown) $\beta$. These graphs also indicate that for  $K_n \ge 4$ the estimate $\widehat\phi_n$ becomes really 
inefficient. Similar facts were observed in the Gaussian case studied in Leipus et al. (2006) and
Celov et al. (2010). Since Figures \ref{fig:miseGam} and \ref{fig:misenotcl} appear rather similar,  
we may conclude that the  differences in the noise distribution and 
the asymptotic results of Section 3 
do not have a strong effect on the performance of the estimators of the mixing density.  
%\end{example}

\begin{figure}
  \centering
 \includegraphics[width=\textwidth,height=6cm]{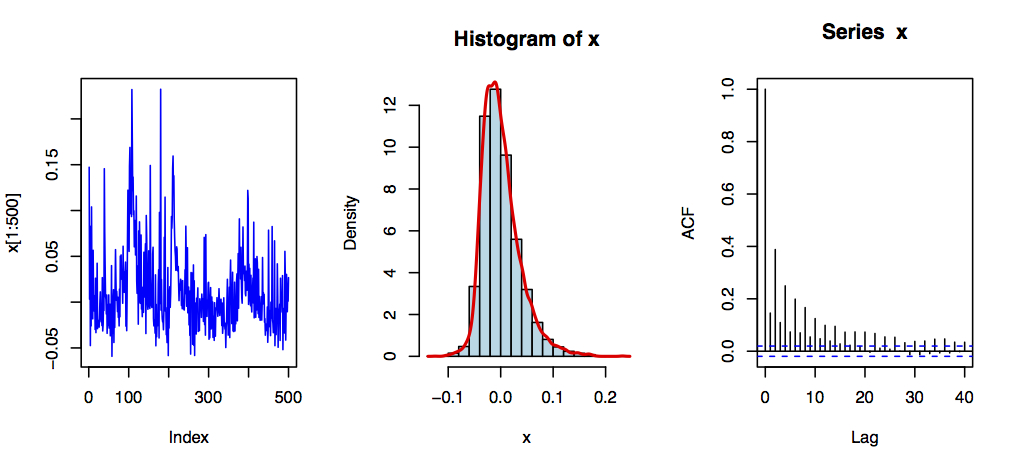}
\caption{The process obtained by aggregating $N=5000$ independent random-coefficient AR(1) with the Gamma noise in
  \eqref{notcl} and mixing density \eqref{exdens}, $\beta = 0.75$. 
%constructed by aggregation of  N=5000 independent processes:
[left] the  first $500$
  values of the simulated
  trajectory,  [Middle] histogram, [right]  empirical auto
  covariance. The sample size  $n=10 000$. }
  %and $\phi(x)
  %\propto \sqrt{x}(1-x)^{0.75} 1_{(0,1)}(x) $.
  %$\beta = 0.75$.}
  \label{fig:traj}
\end{figure}

\begin{figure}
  \centering
 \includegraphics[width=\textwidth,height=8cm]{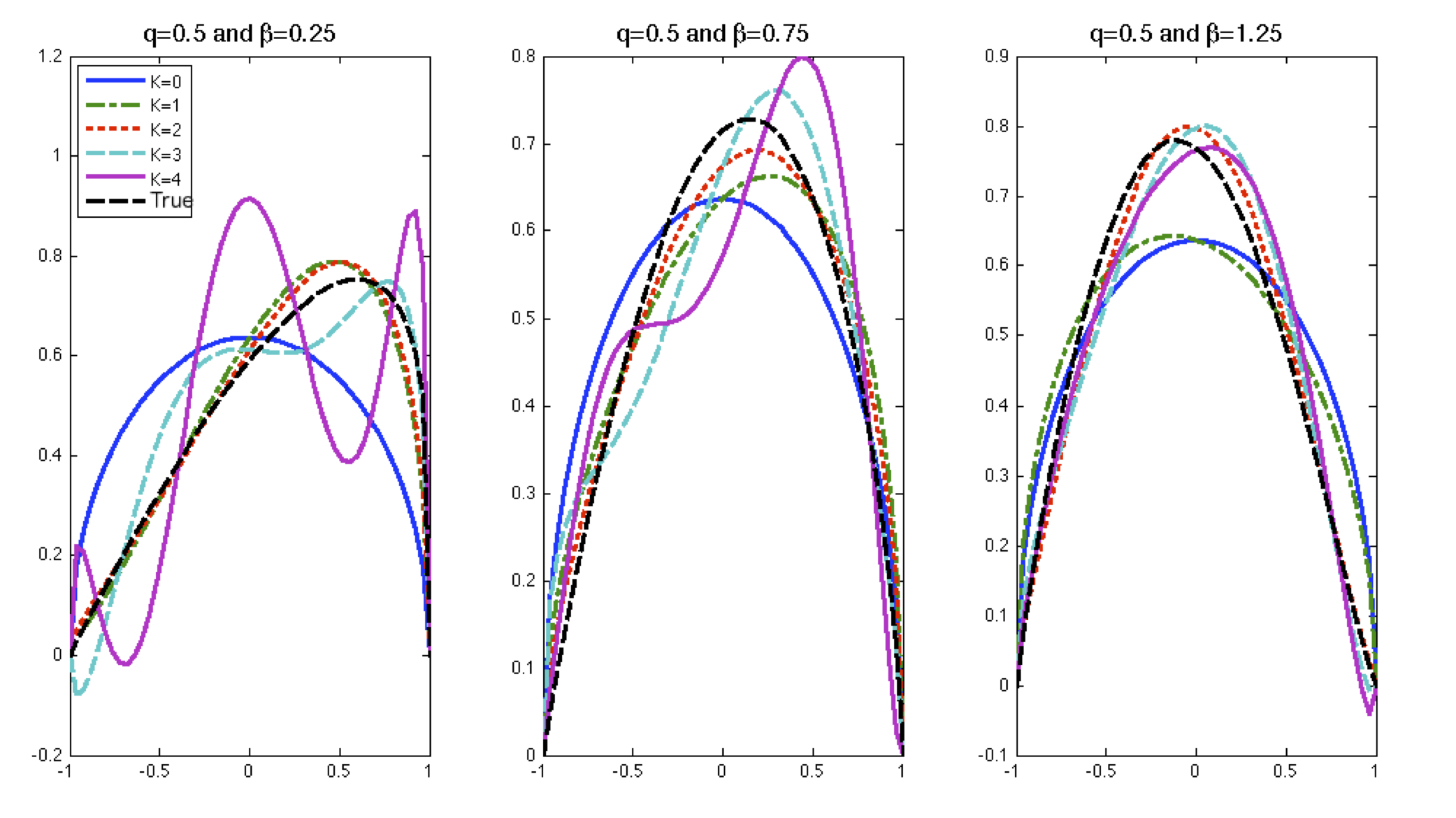}
\caption{The estimates  $\widehat\phi_n$  computed
from the aggregated series with $N=5000$ and Gamma noise 
  \eqref{notcl}. The mixing density is \eqref{exdens}. [left]
  %picture corresponds to a different value of  $\beta$ :
  $\beta = 0.25$,\; [middle]
  $\beta = 0.75$,\; [right] $\beta = 1.25$.
  The sample size $n=10 000$.} %The true density is given by \eqref{exdens}}
  \label{fig:real}
\end{figure}

\begin{figure}
  \centering
\includegraphics[width=\textwidth,height=8cm]{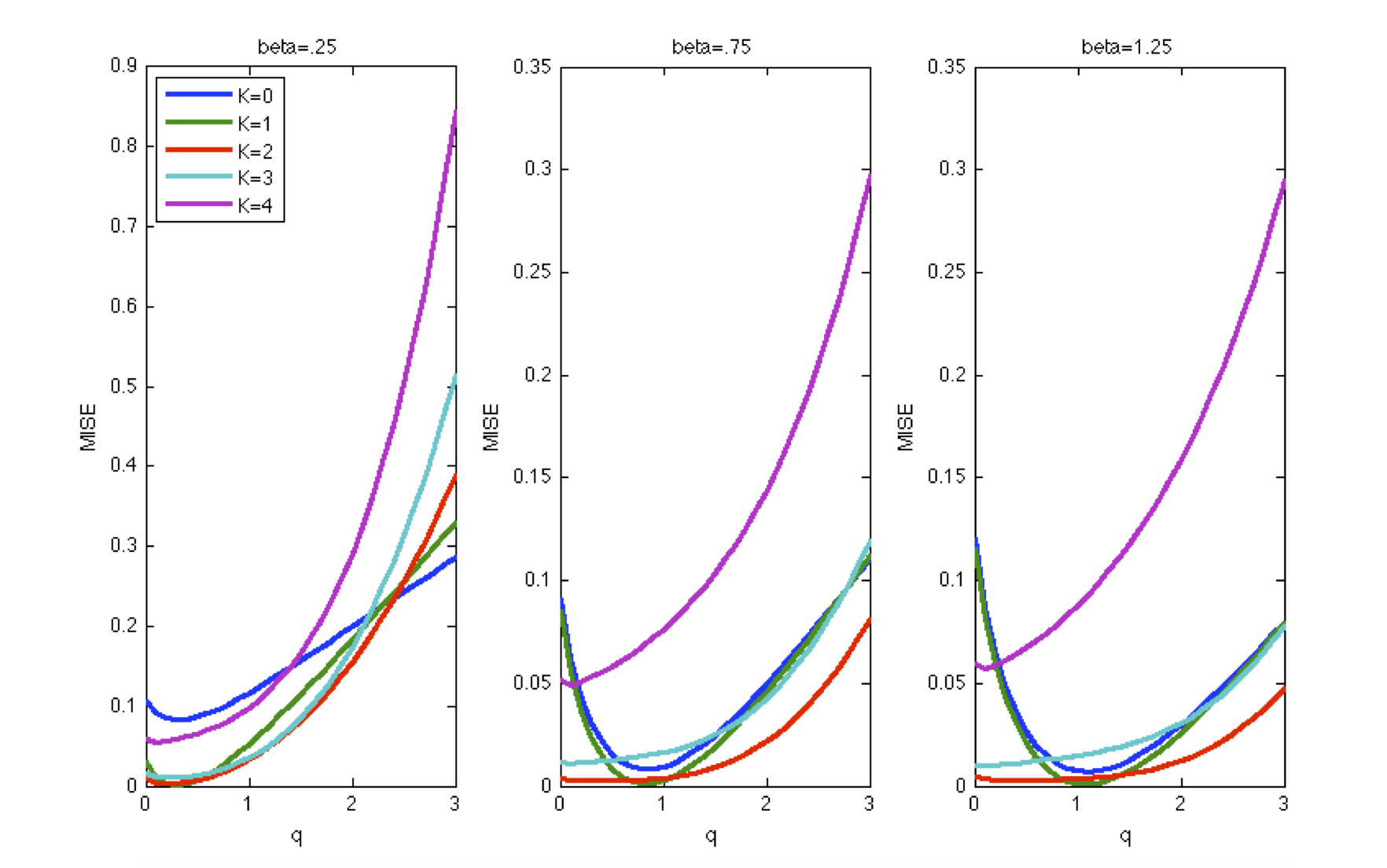}
\caption{The estimated MISE of $\widehat\phi_n$ versus $q$  computed
from the aggregated series with $N=5000$ and the Gamma noise in
\eqref{notcl}. The true density is \eqref{exdens}. [left]
$\beta = 0.25$,\; [middle]
$\beta = 0.75$,\; [right] $\beta = 1.25$. The number of replications 
is 500. The sample size $n=10 000$. }
%Estimation by Monte Carlo method (500 replications) of the
%  MISE as function of $q$ for noise distribution \eqref{notcl}. Each
%  picture corresponds to a different value of  $\beta$ : $\beta = 0.25,\;
%  0.75,\; 1.25$.
%  The sample size of aggregated
%  model is $n=10 000$. The true density is given by \eqref{exdens}}
\label{fig:miseGam}
\end{figure}

\begin{figure}
  \centering
\includegraphics[width=\textwidth,height=8cm]{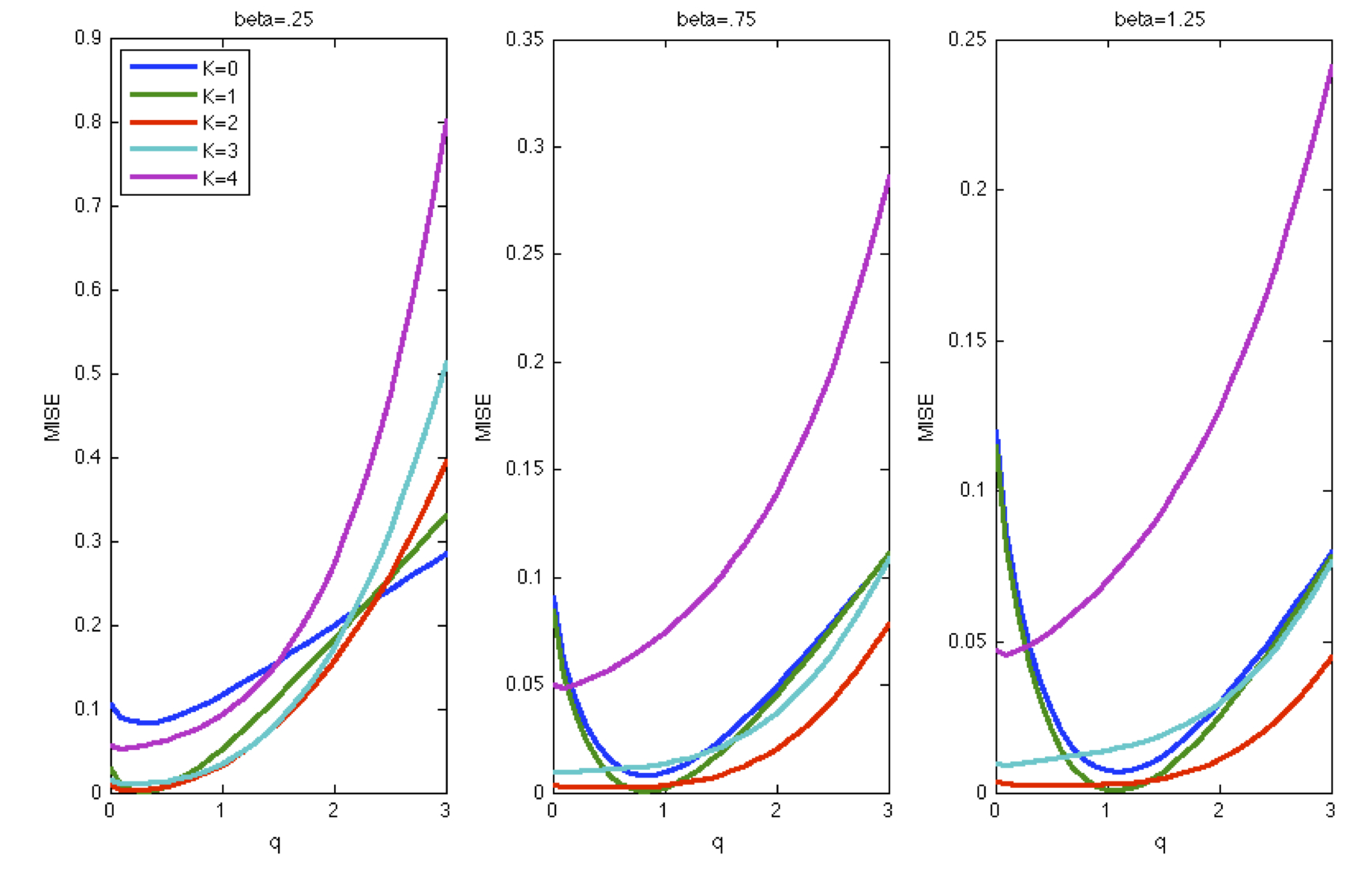}
\caption{The estimated MISE of $\widehat\phi_n$ versus $q$  computed
from the aggregated series with $N=5000$ and  Gaussian noise 
\eqref{gamtcl}. The true density is \eqref{exdens}. [left]
$\beta = 0.25$,\; [middle]
$\beta = 0.75$,\; [right] $\beta = 1.25$. The number of replications
is 500. The sample size $n=10 000$. }
%[Left] Estimation by Monte Carlo method (500 replications) of the
%  MISE as function of $q$ for noise distribution \eqref{gamtcl}. Each
%  picture corresponds to a different value of  $\beta$ : $\beta = 0.25,\;
%  0.75,\; 1.25$.
%  The sample size of aggregated
%  model is $n=10 000$. The true density is given by \eqref{eq:exdens}}
\label{fig:misenotcl}
\end{figure}

\clearpage
\section*{References}
\bigskip

\small

\begin{description}
\itemsep -.04cm

\item{\sc Abramovitz, M. and Stegun, I.}, (1965). {\em Handbook of Mathematical Functions with Formulas, Graphs, and Mathematical Tables.} Dover, New York.

\item {\sc Beran, J., Schuetzner, M. and Ghosh, S.}, (2010). From short to long memory: Aggregation and estimation.
{\em Comput. Stat. Data Anal.} {\bf 54}, 2432--2442.

\item{\sc Bhansali, R.J., Giraitis, L. and Kokoszka, P.}, (2007). Approximations and limit theory for quadratic forms of linear processes. {\em Stoch. Proc. Appl.} {\bf 117}, 71--95.

\item {\sc Billingsley, P.}, (1968). {\em Convergence of Probability Measures.} Wiley, New York.

\item{\sc Brandt, A.}, (1986). The stochastic equation $Y_{n+1}=A_nY_n+B_n$ with stationary coefficients. {\em Adv. Appl. Prob.} {\bf 17}, 211--220.

\item {\sc Celov, D., Leipus, R. and Philippe, A.}, (2007). Time series aggregation, disaggregation and long memory.
{\em Lithuanian Math. J.} {\bf 47}, 379--393.

\item {\sc Celov, D., Leipus, R. and Philippe, A.}, (2010). Asymptotic normality of the mixture density estimator in a disaggregation scheme.
{\em J. Nonparametric Statist.} {\bf 22}, 425--442.

\item{\sc Cox, D.R.}, (1984). Long-range dependence: a review. In: H. A. David and H. T. David (Eds.) {\em Statistics: An Appraisal.} Iowa State University Press, Iowa, 55--74.

\item{\sc Feller, W.}, (1966). {\em An Introduction to Probability Theory and Its Applications}, vol. 2. Wiley, New York.

\item {\sc Gon{\c c}alves, E.  and Gouri{\' e}roux, C.}, (1988). Aggr{\'
e}gation de processus autoregressifs d'ordre 1. {\em Annales
d'Economie et de Statistique} {\bf 12}, 127--149.

\item {\sc Granger, C.W.J.}, (1980). Long memory relationship and the
aggregation of dynamic models. {\em J. Econometrics} {\bf 14},
227--238.

\item{\sc Ibragimov, I.A. and Linnik, Yu.V.}, (1971). {\em  Independent and
Stationary Sequences of Random Variables.} Wolters-Noordhoff,
Groningen.

\item{\sc Lavancier, F.}, (2005). Long memory random fields.
In: P.
Bertail, P. Doukhan, P. Soulier (Eds.), {\em Dependence in
Probability and Statistics.} Lecture Notes in Statistics, vol.
187, pp.195--220. Springer, Berlin.

\item{\sc Lavancier, F.}, (2011). Aggregation of isotropic random fields. {\em J. Statist. Plan. Infer.} {\bf 141}, 3862--3866.

\item{\sc Lavancier, F., Leipus, R. and Surgailis, D.}, (2012). Aggregation of anisotropic random-coefficient autoregressive random field. Preprint.

\item{\sc Leipus, R. and Surgailis, D.} (2003) Random coefficient autoregression, regime
switching and long memory.
{\em Adv. Appl. Probab.} {\bf 35 }, 1--18.

\item{\sc Leipus, R., Oppenheim, G., Philippe, A. and  Viano, M.-C.}, (2006). Orthogonal
series density estimation in a disaggregation scheme.
{\em J. Statist. Plan. Inf.} {\bf 136}, 2547--2571.

\item{\sc Leonenko, N. and Taufer, E.}, (2013). Disaggregation of spatial autoregressive processes.
{\em J. Spatial Statistics} (in press).

\item{\sc Mikosch, T., Resnick, S., Rootz\'en, H. and Stegeman, A.}, (2002). Is network traffic approximated by stable L\'evy motion or fractional Brownian motion? {\em Ann. Appl. Probab.}, {\bf 12}, 23--68.

\item{\sc Puplinskait\.e, D. and Surgailis, D.}, (2012).
Aggregation of autoregressive random fields and anisotropic long memory. Preprint.

\item{\sc Puplinskait\.e, D. and Surgailis, D.}, (2009). Aggregation of random coefficient  AR1(1) process with infinite variance and common innovations. {\em Lithuanian Math. J.} {\bf 49}, 446--463.

\item{\sc Puplinskait\.e, D. and Surgailis, D.}, (2010). Aggregation
of random coefficient  AR1(1) process with infinite variance and idiosyncratic innovations. {\em Adv. Appl. Probab.} {\bf 42}, 509--527.

\item{\sc Reed, M. and Simon, B.}, (1975). {\em Methods of Modern Mathematical Physics}, vol.2. Academic Press, New York.

\item{\sc Robinson, P.} (1978) Statistical inference for a random coefficient
autoregressive model. {\em Scand. J. Statist.} {\bf 5}, 163--168.

\item{\sc  Rajput, B. S. and Rosinski, J.}, (1989). Spectral representations of infinitely
divisible processes. {\em Probab. Th. Rel. Fields} {\bf 82}, 451--487.

 \item {\sc Oppenheim, G. and Viano, M.-C.}, (2004). Aggregation of
random parameters Ornstein-Uhlenbeck or AR processes: some
convergence results. {\em J. Time Ser. Anal.} {\bf 25}, 335--350.

\item{\sc Sato, K.-I.}, (1999). {\em L\'evy processes and Infinitely Divisible Distributions.} Cambridge Univ. Press, Cambridge.

\item{\sc Surgailis, D.} (1981) On infinitely divisible self-similar random fields.
{\em Zeit. Wahrsch.  verw. Geb.}
{\bf 58}, 453--477.

 \item {\sc Zaffaroni, P.} (2004) Contemporaneous aggregation of
linear dynamic models in large economies. {\em J. Econometrics} {\bf
120}, 75--102.

\end{description}

\end{document}